\documentclass[11pt]{article}
\usepackage{latexsym,amssymb,amsmath,amsfonts,amsthm}
\usepackage{graphicx}

\usepackage{amsmath,amssymb}
\usepackage{amsthm,thmtools}
\usepackage{bm}
\usepackage{enumerate}

\usepackage{standalone}
\usepackage{lineno}

\usepackage{tikz}
\newcommand{\R}{{\mathbb R}}
\newcommand{\Z}{{\mathbb Z}}

\newcommand{\C}{{\mathbb C}}

\newcommand{\no}{\nonumber}
\newcommand{\be}{\begin{eqnarray}}
\newcommand{\ben}{\begin{eqnarray*}}
\newcommand{\en}{\end{eqnarray}}
\newcommand{\enn}{\end{eqnarray*}}

\newcommand{\eps}{\epsilon}

\newcommand{\Om}{\Omega}

\newcommand{\al}{\alpha}
\newcommand{\la}{\lambda}

\newtheorem{theorem}{Theorem}[section]
\newtheorem{lemma}[theorem]{Lemma}

\newtheorem{definition}[theorem]{Definition}
\newtheorem{remark}[theorem]{Remark}

\numberwithin{equation}{section} 

\begin{document}
	\renewcommand{\theequation}{\arabic{section}.\arabic{equation}}
	\begin{titlepage}
		\title{\bf Finite element and integral equation methods to conical  diffraction by imperfectly conducting gratings}
		
		\author{Guanghui Hu  \footnote{ School of Mathematical Sciences and LPMC, Nankai  University, Tianjin 300071, People's Republic of China (ghhu@nankai.edu.cn)}\quad 
			\ Jiayi Zhang \footnote{School of Mathematical Sciences and LPMC, Nankai  University, Tianjin 300071, People's Republic of China  (zhangjy97@mail.nankai.edu.cn).}
			\ and \ 	
			Linlin Zhu\footnote{(Corresponding author) School of Mathematical Sciences and LPMC, Nankai  University, Tianjin 300071, People's Republic of China (llzhu@mail.nankai.edu.cn).}			
		}
		\date{}
	\end{titlepage}
	\maketitle
	\vspace{.2in}
	
	\begin{abstract}
		In this paper we study the variational method and integral equation methods for a conical diffraction problem for imperfectly conducting gratings modeled by the impedance boundary value problem of the Helmholtz equation in periodic structures. 
		We justify the strong ellipticity of the sesquilinear form corresponding to the variational formulation and prove the uniqueness of solutions at any frequency. Convergence  of the finite element method using the transparent boundary condition (Dirichlet-to-Neumann mapping) is verified. 
		The boundary integral equation method is also discussed. 
		
		\vspace{.2in} {\bf Keywords:} Diffraction gratings, conical diffraction, variational methods, integral equation methods, finite element analysis, well-posedness.	
	\end{abstract}
	
	\section{Introduction}
	Grating diffraction problems have been extensively studied in the literature via variational and integral equation methods.  We refer to \cite{Bao2022,Petit1980} and references therein
	for mathematical analysis and numerical treatment.  
	In the polarization case, one must assume that the diffraction grating is periodic in one direction ($x_1$-direction), invariant in another direction ($x_3$-direction) and that the incident direction of a time-harmonic electromagnetic plane wave is orthogonal to the $ox_1x_3$-plane.  
	In the TE (resp. TM) polarization, the electric (resp. magnetic) field is parallel to the $x_3$-direction. 
	In this paper we suppose that a time-harmonic plane wave incident obliquely on an imperfectly conducting  grating, which leads to the so-called conical diffraction problems.
	The impedance boundary condition will be used to model the boundary behavior of the wave fields between a highly conducting material and an isotropic, homogeneous and lossless background medium.  
	
	In periodic structures,	Elschner, Hinder, Penzel and Schmidit \cite{El00} proved the well-posedness of the conical diffraction problem with transmission conditions via the variational method. 
	Elschner and Schmidt \cite{Elschner2003} studied stability of the conical diffraction problem with respect to variation of the grating profile and obtained explicit formulas for the derivatives of reflection and transmission coefficients with respect to perturbations of interfaces.	
	If the scattering object is an infinitely long cylinder, the conical diffraction is also referred to as oblique scattering.
	In \cite{Nakamura2013}, Wang and Nakamura applied the integral equation method to prove the well-posedness of the oblique scattering problem in a homogeneous medium. 
	In an inhomogeneous medium, the uniqueness and existence of the oblique problem were also studied through the Lax-Phillips method; see Nakamura \cite{HaibingWang2012}. 
	In this paper, we consider the conical diffraction problem in periodic structures under the impedance boundary condition and investigate both finite element and boundary integral equation methods.
	
	The outline of the paper is organized as follows. 
	In Section \ref{sec:conical}, we formulate the conical diffraction problem by deriving a coupled Helmholtz system with the impedance boundary condition from Maxwell's system. 
	In Section \ref{sec:VF}, we state the variational formulation in one periodic cell with the DtN operator imposed on the artificial boundary. 
	An energy formula is verified to prove  the uniqueness of the truncated boundary value problem. 
	The strong ellipticity of the variational formulation is shown and the well-posedness of the diffraction problem follows from the Fredholm theory. 
	In Section \ref{sec:FEM}, we show the convergence of the finite element method based on the variational formulation. 
	Finally, the integral equation method will be briefly discussed in Section \ref{sec:integralEm}.

	\section{Conical diffraction problem} \label{sec:conical}
\begin{figure}[htp]
		\centering
		\begin{tikzpicture}
	\tikzstyle{surffill} = [fill=blue!20,fill opacity=0.8];

	\coordinate (O) at (0,0,0);
	\draw[thick,dashed] (O) -- (6,0,0);
	\draw[thick,dashed] (O) -- (0,2,0);
	\draw[thick,dashed] (O) -- (0,0,2);
	\draw[thick,->] (6,0,0) -- (8,0,0) node [right] {$x_{1}$};
	\draw[thick,->] (0,2,0) -- (0,4,0) node [left] {$x_{2}$};
	\draw[thick,->] (0,0,2) -- (0,0,4) node [left] {$x_{3}$};
	
	\coordinate (1) at (3,3,0);
	\draw[thick,->] (1) -- (4,3,0) node [right] {\tiny$x_{1}$};
	\draw[thick,->] (1) -- (3,4.25,0) node [left]  {\tiny$x_{2}$};
	\draw[thick,->] (1) -- (3,3,1) node [left]  {\tiny$x_{3}$};

	\draw[-stealth,cyan] (2,5,2) -- (2.5,4,1) node [below left] {$u^{i}$};
	\draw[very thin,gray] (2.45,4.1,0) -- (1);

	\draw[fill=green!20,ultra thin,opacity=0.6,dashed] (1) -- (2.5,4,1) -- (2.5,4,0) -- cycle;
	\draw[fill=pink!20,ultra thin,opacity=0.6,dashed] (1) -- (2.5,4,0) -- (3,4,0) -- cycle;

	\draw[surffill] (0,2,0) .. controls (2,3,0) and (4,1,0) .. (6,2,0) 
	                -- (6,2,2) .. controls (4,1,2) and (2,3,2) .. (0,2,2) 
	                -- cycle;
	\draw[surffill] (0,2,2) .. controls (2,3,2) and (4,1,2) .. (6,2,2) 
		             -- (6,0,2) -- (0,0,2) -- cycle;
	\draw[surffill] (6,2,0) -- (6,2,2) -- (6,0,2) -- (6,0,0) -- cycle;
	
	\draw       (5,2.5,0)    node [above] {$\tilde{\Omega}$};
	\draw       (4,2,1.5)    node [right] {$\tilde{\Gamma}$};
	\draw       (1,2.5,0)    node [above] {$\epsilon$,$\mu$};
	\draw       (6,1,1)      node [right] {\Large$\ldots$};
	\draw       (0,1,1)      node [left]  {\Large$\ldots$};
	\draw       (3.1,3.5,0)  node [left]  {\scriptsize$\theta$};
	\draw       (2.9,3.35,0) node [left]  {\scriptsize$\phi$};
\end{tikzpicture}
		\caption{Geometry of the three-dimensional conical diffraction problem in one periodic cell. $ \phi $ is the angle between incident direction $ \textbf{k} $ and $(x_1,x_2)$-plane. $ \theta $ is the angle between $ (\alpha,-\beta) $ and the $x_2$-axis.}
		\label{fig:geo3d}
	\end{figure}
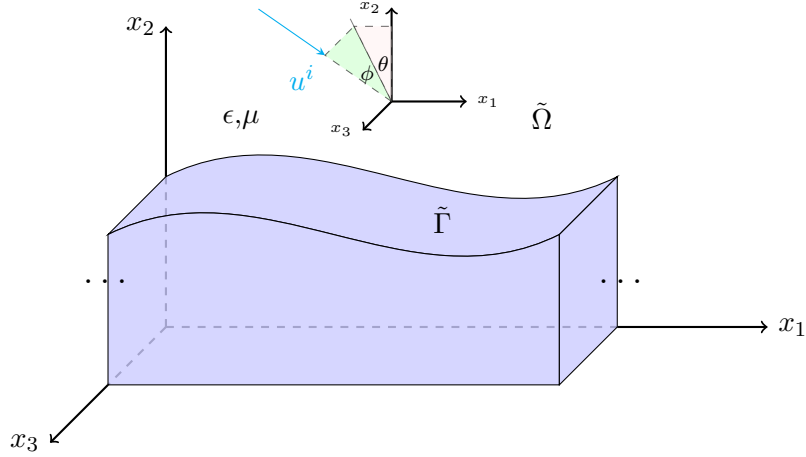

	Assume an incoming time-harmonic plane wave of the form
	\be\label{incoming}
	(\mathcal{E}^{in},\mathcal{H}^{in})=(\textbf{p}e^{i\al x_1-i\beta x_2+i\gamma x_3},
	\textbf{q}e^{i\al x_1-i\beta x_2+i\gamma x_3}) e^{-i\omega t}=:(\textbf{E}^{in},\textbf{H}^{in})e^{-i\omega t},
	\en
	is incident on an imperfectly conducting grating with a high conductivity embedded in an isotropic homogeneous medium in $\R^3$. 
	Denote by $\tilde{\Gamma}$ the grating profile and $\tilde{\Om}$ the unbounded domain above $\tilde{\Gamma}$. 
	The diffraction problem  can be modeled by the reduced Maxwell's system
	\be\label{Maxwell}
	\nabla\times \textbf{E}=i\omega \mu \textbf{H}, \quad \nabla\times\textbf{H}=-i\omega  \eps \textbf{E}\quad\mbox{in}\quad\tilde{\Omega} ,
	\en
	where the total fields $(\textbf{E},\textbf{H})$ are the sum of the incident waves $(\textbf{E}^{in}, \textbf{H}^{in})$ and the outgoing scattered waves $(\textbf{E}^{sc}, \textbf{H}^{sc})$ in $\tilde{\Om}$.
	In (\ref{Maxwell}), $\omega$ denotes the angular frequency. 
	The dielectric coefficient $\eps$ and the magnetic permeability $\mu$ of  the homogeneous medium in $\tilde{\Omega}$ are both assumed to be positive constants. 
	Set $k = \omega\sqrt{\epsilon\mu}$ as the wavenumber of the background medium.
	We enforce the impedance boundary condition on $ \tilde{\Gamma} $
	\be\label{impedance}
	\nu\times \textbf{E}\times\nu=\lambda\; (\nu\times \textbf{H})\quad\mbox{on}\quad \tilde{\Gamma}, 
	\en 
	where $\nu=(\nu_1,\nu_2,\nu_3)\in \mathbb{S}^2:=\{x\in\R^3: |x|=1\}$ is normal to $\tilde{\Gamma}$ directed into the exterior of $\tilde{\Omega}$ and $\lambda<0$ is the impedance coefficient which is assumed to be a constant. 
	The problem (\ref{incoming})-(\ref{impedance}) is called a conical diffraction problem if the incident direction $ \textbf{k}=:(\alpha,-\beta,\gamma) $ is not orthogonal to the $x_{2}$-direction, i.e., $\gamma \neq 0$.
	For conical diffraction problems, the wave vectors of the reflected or transmitted propagating modes lie on the surface of a cone whose axis is parallel to the $x_{3}$-direction \cite{Schmidt11}. 
	We refer to Figure \ref{fig:geo3d} for an illustration of the grating conical diffraction problem. 
	
	In order for $(\textbf{E}^{in}, \textbf{H}^{in})$ given in (\ref{incoming}) to satisfy (\ref{Maxwell}), the constant amplitude vector $\textbf{p}$ must be perpendicular to the wave vector $\textbf{k}=(\alpha,-\beta,\gamma)$, that is $\textbf{p}\cdot \textbf{k}=0$. 
	Furthermore $\textbf{k}\cdot \textbf{k}=k^2=\omega^2\eps\mu$ and $\textbf{q}=(\omega\mu)^{-1} \textbf{k}\times \textbf{p}$. 
	We can express the wave vector $\textbf{k}$ as \ben
	\textbf{k}=(\alpha,-\beta,\gamma):=k(\sin\theta\cos\phi, -\cos\theta\cos\phi,\sin\phi),
	\enn
	in terms of the angles of incidence $\theta,\phi\in(-\pi/2,\pi/2)$.
	
\begin{figure}[htp]
		\centering
		\begin{tikzpicture}
	\tikzstyle{surffill} = [fill=blue!20,fill opacity=0.8];
	
	\draw[dashed] (0,2) -- (0,4);
	\draw[dashed] (6,2) -- (6,4);
	
	\draw[surffill] (0,2) .. controls (2,3) and (4,1) .. (6,2) 
	                -- (6,0) -- (0,0) -- cycle;
	\draw (4,2.5) node [above] {$\Omega$};
	\draw (4,2)   node [right] {$\Gamma$};
	\draw (2,2.5) node [above] {$\epsilon$, $\mu$};	
	\draw (0,4) node [left]  {$x_{2}$};
	\draw (6,0) node [right] {$x_{1}$};
\end{tikzpicture}
		\caption{Geometry of the conical diffraction problem.}
		\label{fig:geo2d}
	\end{figure}
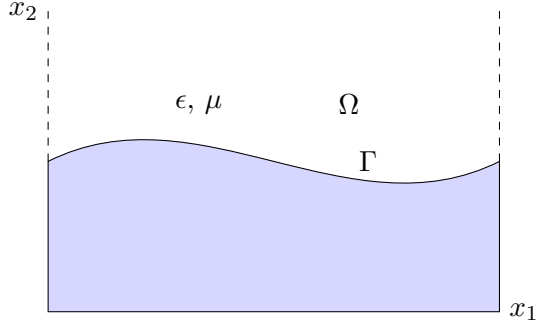

	Assume that $\tilde{\Gamma}$ remains invariant in $x_3$ and is $2\pi$-periodic in $x_1$. 
	If the incoming wave is of the form (\ref{incoming}), we make an ansatz on the total field
	\ben
	(\textbf{E},\textbf{H})(x_1,x_2,x_3)=(E(x_1,x_2),H(x_1,x_3))\;e^{i\gamma x_3},  
	\enn 
	with $E=(E_1,E_2,E_3), H=(H_1,H_2,H_3): \R^2\rightarrow\C^3$. 
	The Maxwell equations (\ref{Maxwell}) can be reduced to two Helmholtz equations for the total fields $u=E_3$ and $v=H_3$ (see \cite{El00}):
	\ben
	\Delta u+\kappa^2u=0,\quad \Delta v+\kappa^2 v=0&&\mbox{in}\quad\Omega,\quad \kappa^2=k^2-\gamma^2. 
	\enn
	Here $\Omega$ denotes the restriction of the cross-section of $\tilde{\Omega}$ by the $(x_1,x_2)$-plane to one periodic cell $(0,2\pi)$. 
	Analogously, denote by $\Gamma$ the  counter part of $\tilde{\Gamma}$ in the periodic cell $(0,2\pi)$.
	The reduced geometry of the conical diffraction problem is shown in Figure \ref{fig:geo2d}.
	Next, we turn to the reduction of the boundary condition (\ref{impedance}) in $\R^2$. 
	Obviously, we have $\nu_3=0$ and
	\be\label{I4}
	(\nu\times E)\times \nu&=&(-\nu_2(\nu_1E_2-\nu_2E_1),\,-\nu_1(\nu_1E_2-\nu_2E_1),\, E_3),\\ \no
	\nu\times H&=&(\nu_2H_3,\;-\nu_1H_3,\;\nu_1H_2-\nu_2H_1).
	\en
	Moreover, there holds (see \cite{El00})
	\be\label{I3}
	\nu_1E_2-\nu_2E_1=\frac{i\gamma}{\kappa^2}\frac{\partial E_3}{\partial \tau}-\frac{i\omega\mu}{\kappa^2}\frac{\partial H_3}{\partial n},\quad
	\nu_1H_2-\nu_2H_1=\frac{i\gamma}{\kappa^2}\frac{\partial H_3}{\partial \tau}+\frac{i\omega\eps}{\kappa^2}\frac{\partial E_3}{\partial n},
	\en 
	with 
	$$n=(\nu_1,\nu_2)\in\mathbb{S},\;\tau=(-\nu_2,\nu_1)\in\mathbb{S},\quad\partial_n=\nu_1\partial_{1}+\nu_2\partial_2,\; \partial_\tau=-\nu_2\partial_{1}+\nu_1\partial_2,\;\partial_j=\frac{\partial}{\partial x_j}. $$
	Meanwhile, for the reduced Helmholtz equation, the incoming time-harmonic plane wave (\ref{incoming}) takes the form
	\begin{equation*}
		u^i=p_3 e^{i\alpha x_1+i\beta x_2},\quad v^i=p_3 e^{i\alpha x_1+i\beta x_2}.
	\end{equation*}
	Combining (\ref{I4})-(\ref{I3}) and the impedance boundary condition (\ref{impedance}), we get
	\ben
	-\frac{i\gamma}{\kappa^2}\frac{\partial E_3}{\partial \tau}+\frac{i\omega\mu}{\kappa^2}\frac{\partial H_3}{\partial n}=\lambda H_3,\quad
	E_3=\lambda\,\left(\frac{i\gamma}{\kappa^2}\frac{\partial H_3}{\partial \tau}+\frac{i\omega\eps}{\kappa^2}\frac{\partial E_3}{\partial n}\right),
	\enn
	which, for $u=E_3$ and $v=H_3$, is equivalent to the boundary condition
	\be\label{eqn:BC}
	\la\frac{\partial u}{\partial n}+\frac{i\kappa^2}{\omega\eps}u+\frac{\la\gamma}{\omega\eps}\frac{\partial v}{\partial\tau} =0,\quad
	\frac{\partial v}{\partial n}+\frac{i\la \kappa^2}{\omega\mu}v-\frac{\gamma}{\omega\mu}\frac{\partial u}{\partial \tau}=0\quad\mbox{on}\quad\Gamma. 
	\en
	Using $\gamma=\omega\sqrt{\eps\mu}\sin\phi$ and $\kappa^2=k^2\cos^2\phi=\omega^2\mu\eps \cos^2\phi$, the previous boundary condition can be written as
	\be\label{BC1}\left\{\begin{array}{lll}
		\la\frac{\partial u}{\partial n}+i\omega\mu \cos^2\phi\; u+\la\sin\phi\sqrt{\frac{\mu}{\eps}}\frac{\partial v}{\partial\tau} &=&0,\\
		\frac{\partial v}{\partial n}+i\la\omega\eps\cos^2\phi\, v-\sin\phi\sqrt{\frac{\eps}{\mu}}\frac{\partial u}{\partial \tau}&=&0,\end{array}\right. \quad\mbox{on}\quad\Gamma.
	\en
	\begin{remark}
		If $\la=0$, then the boundary condition (\ref{eqn:BC}) (or (\ref{BC1})) reduces to $\frac{\partial v}{\partial n}=u=0$, which corresponds to the TE or TM polarization of the electromagnetic scattering by perfectly conducting gratings. 
		If $\phi=0$, then both $u$ and $v$ satisfy the standard impedance boundary condition for the Helmholtz equation.
	\end{remark}
	
	\section{Radiation condition and variational formulation} \label{sec:VF}
	For $b>\Gamma_{\max}:=\max_{x\in \Gamma}\{x_2\}$, define
	\ben
	\Gamma_b:=\{(x_1,b): 0<x_1<2\pi\},\quad \Omega_b:=\{x\in\Omega: x_2<b\}. 
	\enn
	A function $ u(x_1,x_2) $ is called $ \alpha $-quasiperiodic if $ e^{-i\alpha x_1}u(x_1,x_2) $ is $ 2\pi $-periodic in $ x_1 $, or equivalently, 
	\ben
	u(x_1+2\pi,x_2)=e^{2i\alpha \pi}u(x_1,x_2).
	\enn
	Since the incident field is $ \alpha  $-quasiperiodic, the scattered field $ u^s,\; v^s $ are also assumed to be $ \alpha $-quasiperiodic. 
	Then the function $  u^s(x_1,x_2)e^{-i\alpha x_1},\; v^s(x_1,x_2)e^{-i\alpha x_1}$ can be expended as a Fourier series.
	Inserting these series into the Helmholtz equation, we can express $ u^s $ and $ v^s $ as a sum of plane waves. 
	Physically, the scattered field $\left(u^{s}, v^{s}\right)$ remain bounded as $x_2 \rightarrow \infty$, leading to the well-known Rayleigh expansion condition:
	\begin{align}
		u^s(x)=\sum_{n\in \mathbb{Z}}u_n e^{i\alpha_n x_1+i\beta_n x_2},\qquad 
		v^s(x)=\sum_{n\in \mathbb{Z}}v_n e^{i\alpha_n x_1+i\beta_n x_2},\quad x_2>\Gamma_{\max}, \label{eqn:UVs}
	\end{align}
	with the Rayleigh coefficients $u_n,v_n\in \mathbb{C}$, where
	\begin{align*}
		\alpha_n:=n+\alpha,  \qquad
		\beta_n:=\begin{cases}
			\sqrt{\kappa^2-|\alpha_n|^2},& |\alpha_n|\leq\kappa,\\
			i\sqrt{|\alpha_n|^2-\kappa^2}, &  |\alpha_n|> \kappa,
		\end{cases}
	\end{align*}
	with $i=\sqrt{-1}$.
	It is clear that $(u^s, v^s)$ in (\ref{eqn:UVs}) can be split into the finite sum $\sum_{|\alpha_n|\leq k}$ of outgoing plane waves and the infinite sum $\sum_{|\alpha_n|> k}$ of exponentially decaying waves, which are called surface or evanescent waves. We summarize our conical diffraction problem as follows:
	\be\label{cdp}  
	\qquad\left\{\begin{array}{ll}
		\Delta u+\kappa^2u=0,\ \Delta v+\kappa^2 v=0,  &\mbox{in} \quad \Omega, \\
		\la\frac{\partial u}{\partial n}+\frac{i\kappa^2}{\omega\eps}u+\frac{\la\gamma}{\omega\eps}\frac{\partial v}{\partial\tau} =0,\ \frac{\partial v}{\partial n}+\frac{i\la \kappa^2}{\omega\mu}v-\frac{\gamma}{\omega\mu}\frac{\partial u}{\partial \tau}=0,     &\mbox{on} \quad \Gamma,\\
		\mbox{$u^s$ and $v^s$ fulfill the Rayleigh expansion \eqref{eqn:UVs}.}				
	\end{array}\right.  
	\en
	Then we introduce the variational space
	\ben
	X=\{(u,v)\in H^1(\Om_b)^2: u,v\; \mbox{are $\al$-quasiperiodic}\}.
	\enn
	In order to derive the variational formula, we will need Green's formula for functions in $ H_\alpha^1(\Omega_b) $, for which it is well known.
	\begin{lemma}
		Assume that $ f\in H_\alpha^2(\Omega_b) $ and $ g\in H_\alpha^{1}(\Omega_b) $, Then
		\ben
		\int_{\Omega_b}\nabla f\cdot \nabla \overline{g}+\Delta f\bar{g}\,dx=\int_{\partial \Omega_b}\partial_n f \overline{g}\,ds,\quad 
		\int_{\Omega_b}\nabla f\cdot \nabla ^\bot\overline{g}\,dx=-\int_{\partial \Omega_b}\partial_\tau f \overline{g}\,ds,
		\enn
		where $\nabla=(\partial_1,\partial_2)$ and $\nabla^\bot=(-\partial_2,\partial_1)$.
	\end{lemma}
	Let $ u,v\in H_\alpha^1(\Omega_b) $ solve the conical diffraction problem (\ref{cdp}). 
	Applying Green's formula to the Helmholtz equations yields
	\begin{align}
		0=\int_{\Omega_b}(\Delta u+\kappa^2u)\overline{\varphi}\,dx&=\int_{\Omega_b}-\nabla u\cdot\nabla\overline{\varphi}+\kappa^2u\bar{\varphi}\,dx+\int_{\partial\Omega_b}\partial_n u\,\overline{\varphi}\,ds,\label{VF1}\\
		\int_{\Omega_b}\nabla v\cdot \nabla^\bot\overline{\varphi}\,dx&=-\int_{\partial\Omega_b}\partial_\tau v\,\overline{\varphi}\,ds \quad\text{ for all } \varphi \in H_\alpha^1(\Omega_b).\label{VF2}
	\end{align}
	Multiplying the equations $ (\ref{VF1}) $ and $ (\ref{VF2}) $ by the constant factors $ \frac{\omega\eps}{\kappa^2} $ and $ \frac{\gamma}{\kappa^2} $, respectively, and taking the difference of the resulting formulas, we get
	\begin{equation}
		\int_{\partial\Omega_b}\frac{\omega\eps}{\kappa^2}\partial_n u\,\overline{\varphi}+\frac{\gamma}{\kappa^2}\partial_\tau v\,\overline{\varphi}\,ds=\int_{\Omega_b}\left[\frac{\omega\eps}{\kappa^2}\nabla u\cdot \nabla\overline{\varphi}-\frac{\gamma}{\kappa^2}\nabla v\cdot \nabla^\bot\overline{\varphi}-\omega \eps\,u\;\overline{\varphi}\right]\;dx. \label{proVF1}
	\end{equation}
	Similarly, we get
	\begin{align}
		0=\int_{\Omega_b}(\Delta v+\kappa^2v)\overline{\psi}\,dx&=\int_{\Omega_b}-\nabla v\cdot\nabla\overline{\psi}+\kappa^2v\bar{\psi}\,dx+\int_{\partial\Omega_b}\partial_n v\,\overline{\psi}\,ds,\label{VF3}\\
		\int_{\Omega_b}\nabla u\cdot \nabla^\bot\overline{\psi}\,dx&=-\int_{\partial\Omega_b}\partial_\tau u\,\overline{\psi}\,ds, \quad\text{ for all } \psi \in H_\alpha^1(\Omega_b). \label{VF4}
	\end{align}
	Multiplying the equations $ (\ref{VF3}) $ and $ (\ref{VF4}) $ by the constant factors $ \frac{\omega\mu}{\kappa^2} $ and $ \frac{\gamma}{\kappa^2} $, respectively, then taking the sum of the two formulas, we get
	\begin{align}
		\int_{\partial\Omega_b}\frac{\omega\mu}{\kappa^2}\partial_n v\,\overline{\psi}-\frac{\gamma}{\kappa^2}\partial_\tau u\,\overline{\psi}\,ds=\int_{\Omega_b}\left[\frac{\omega\mu}{\kappa^2}\nabla v\cdot \nabla\overline{\psi}+\frac{\gamma}{\kappa^2}\nabla u\cdot \nabla^\bot\overline{\psi}-\omega \mu\,v\;\overline{\psi}\right]\;dx.\label{proVF2}
	\end{align}
	Recalling the boundary conditions (\ref{eqn:BC}) on $\Gamma$, the left-hand terms of (\ref{proVF1}) and (\ref{proVF2}) over the integral $\Gamma$ can be reformulated as
	\begin{align*}
		\int_{\Gamma}\frac{\omega\eps}{\kappa^2}\partial_n u\,\overline{\varphi}+\frac{\gamma}{\kappa^2}\partial_\tau v\,\overline{\varphi}\,ds=\int_{\Gamma}\frac{\omega\eps}{\lambda\kappa^2}\left( -\frac{i\kappa^2}{\omega\eps}u\right)\overline{\varphi}\,ds= -\frac{i}{\lambda}\int_{\Gamma} u\,\overline{\varphi}\,ds,\\
		\int_{\Gamma}\frac{\omega\mu}{\kappa^2}\partial_n v\,\overline{\psi}-\frac{\gamma}{\kappa^2}\partial_\tau u\,\overline{\psi}\,ds=\int_{\Gamma}\frac{\omega\mu}{\kappa^2}\left( -\frac{i\la \kappa^2}{\omega\mu}v\right) \overline{\psi}\,ds=-i\lambda \int_{\Gamma} v\,\overline{\psi}\,ds.
	\end{align*}
	Therefore, we need to find $(u,v)\in X$ such that for all $ (\varphi, \psi) \in X $,
	\begin{align}
		0&=\frac{i}{\lambda}\int_{\Gamma} u\,\overline{\varphi}\,ds +\int_{\Omega_b}\left[\frac{\omega\eps}{\kappa^2}\nabla u\cdot \nabla\overline{\varphi}-\frac{\gamma}{\kappa^2}\nabla v\cdot \nabla^\bot\overline{\varphi}-\omega \eps\,u\;\overline{\varphi}\right]\;dx \nonumber \\
		&\quad - \int_{\Gamma_b}\left[\frac{\omega\eps}{\kappa^2}\partial_n u\,\overline{\varphi}+\frac{\gamma}{\kappa^2}\partial_\tau v\,\overline{\varphi}\right]\,ds, \label{v1}\\ 
		0&=i\lambda \int_{\Gamma} v\,\overline{\psi}\,ds +\int_{\Omega_b}\left[\frac{\omega\mu}{\kappa^2}\nabla v\cdot \nabla\overline{\psi}+\frac{\gamma}{\kappa^2}\nabla u\cdot \nabla^\bot\overline{\psi}-\omega \mu\,v\;\overline{\psi}\right]\;dx \nonumber \\  
		&\quad - \int_{\Gamma_b}\left[\frac{\omega\mu}{\kappa^2}\partial_n v\,\overline{\psi}-\frac{\gamma}{\kappa^2}\partial_\tau u\,\overline{\psi}\right]\,ds.\label{v2}
	\end{align} 
	Combining (\ref{v1}) and (\ref{v2}), we get
	\begin{align}
		\int_{\Gamma}\frac{i}{\lambda}u\,\overline{\varphi}+i\lambda v\,\overline{\psi}\,ds+\int_{\Omega_b}\left[\frac{\omega\eps}{\kappa^2}\nabla u\cdot \nabla\overline{\varphi}-\frac{\gamma}{\kappa^2}\nabla v\cdot \nabla^\bot\overline{\varphi}-\omega \eps\,u\;\overline{\varphi}+\frac{\omega\mu}{\kappa^2}\nabla v\cdot \nabla\overline{\psi}\right.\nonumber\\
		+\left.\frac{\gamma}{\kappa^2}\nabla u\cdot \nabla^\bot\overline{\psi}-\omega \mu\,v\;\overline{\psi}\right]\;dx
		-\int_{\Gamma_b}\frac{1}{\kappa^2}
		\left(\begin{array}{ccc}   
			\omega\eps\partial_n u+ \gamma \partial_\tau v \\  
			\omega\mu\partial_n v- \gamma \partial_\tau u \\  
		\end{array}
		\right)
		\cdot
		\left( \begin{array}{ccc}
			\overline{\varphi}\\
			\overline{\psi}
		\end{array}\right)\,ds=0.	\label{VFPre}
	\end{align}	
	\begin{definition}[DtN map] \label{dtn}
		The Dirichlet-to-Neumann (DtN) map $ T $ is defined by
		$$ T:(g_1,g_2)^\top\to -\left(\frac{\omega\eps}{\kappa^2}\partial_n w_1+\frac{\gamma}{\kappa^2}\partial_\tau w_2,
		\frac{\omega\mu}{\kappa^2}\partial_n w_2-\frac{\gamma}{\kappa^2}\partial_\tau w_1\right)^\top\quad \text{on} \;\Gamma_b,$$
		where $w_j (j = 1, 2)$ is the unique radiation solution to the Helmholtz equation $ \Delta w_j +\kappa^2w_j=0 $ in $x_2 > b$ with the Dirichlet boundary condition $w_j=g_j$ on $\Gamma_b$.
	\end{definition}
	Now we want to derive an analytical expression of the DTN map $T$. 
	For the $ \alpha $-quasiperiodic vector function $ g=(g_1,g_2)^\top \in H_\alpha^{1/2}(\Gamma_b)^2$, we can get its Fourier expansion $ g(x_1)=\sum_{n\in\mathbb{Z}}\hat{g}_n e^{i\alpha_n x_1} $, where $\hat{g}_n=(\hat{g}_{n,1},\hat{g}_{n,2})^\top$. 
	It is easy to deduce that
	\ben
	w_j(x)=\sum_{n\in\mathbb{Z}}\hat{g}_{n,j} e^{i\alpha_n x_1+i\beta_n (x_2-b)},\quad  x_2 > b,\ j = 1, 2,
	\enn
	where $ w_j $ is the function specified in the Definition \ref{dtn}.
	Direct calculations show
	\ben
	&\,&-\left(\frac{\omega\eps}{\kappa^2}\partial_n w_1+\frac{\gamma}{\kappa^2}\partial_\tau w_2,
	\frac{\omega\mu}{\kappa^2}\partial_n w_2-\frac{\gamma}{\kappa^2}\partial_\tau w_1\right)\bigg|_{\Gamma_b}\\
	&=&-\frac{1}{\kappa^2}\left(\omega\eps\sum_{n\in\mathbb{Z}}i\beta_n\hat{g}_{n,1} e^{i\alpha_n x_1}+\gamma\sum_{n\in\mathbb{Z}}(-i\alpha_n )\hat{g}_{n,2} e^{i\alpha_n x_1},\right.\\
	&\qquad&\left.\omega\mu\sum_{n\in\mathbb{Z}}i\beta_n\hat{g}_{n,2} e^{i\alpha_n x_1}-\gamma\sum_{n\in\mathbb{Z}}(-i\alpha_n )\hat{g}_{n,1} e^{i\alpha_n x_1}\right)\\
	&=&-\frac{1}{\kappa^2}\sum_{n\in\mathbb{Z}}\left(                 
	\begin{array}{ccc}   
		i\omega\eps\beta_n & -i\gamma \alpha_n \\  
		i\gamma \alpha_n & i\omega\mu\beta_n \\  
	\end{array}
	\right)
	\left( \begin{array}{ccc}
		\hat{g}_{n,1}\\
		\hat{g}_{n,2}
	\end{array}\right)e^{i\alpha_n x_1}\\
	&=&\sum_{n\in\mathbb{Z}}M_n\hat{g}_n e^{i\alpha_n x_1},
	\enn
	where 
	\begin{equation}       
		M_n=\frac{1}{\kappa^2}\left(                 
		\begin{array}{ccc}   
			-i\omega\eps\beta_n & i\gamma \alpha_n \\  
			-i\gamma \alpha_n & -i\omega\mu\beta_n \\  
		\end{array}
		\right).  \label{Mn}         
	\end{equation}
	
	Hence, the operator $ T $ acting on the $\alpha$-quasiperiodic vector function $w \in H_{\alpha}^{1/2}\left(\Gamma_{b}\right)^2$ can be expressed as
	\begin{align*}
		(Tw)(x)=\sum_{n\in\mathbb{Z}}M_n\hat{w}_n e^{i\alpha_n x},   \quad   \hat{w}_n=\frac{1}{2\pi}\int_{0}^{2\pi}w(x)e^{-i\alpha_n x}\,dx \in \mathbb{ C} ^2.
	\end{align*}
	
	\begin{lemma} (see \cite{Bao1995})\label{ctnT}
		The DtN operator $ T:H_\alpha^{1/2}(\Gamma_b)^2 \to H_\alpha^{-1/2}(\Gamma_b)^2 $ is continuous, i.e., there exists a positive constant $ C $ such that
		\ben
		\|Tw\|_{ H_\alpha^{-1/2}(\Gamma_b)^2}\leq C \|w\|_{ H_\alpha^{1/2}(\Gamma_b)^2} \qquad \text{for all } \,w \in H_\alpha^{1/2}(\Gamma_b)^2.
		\enn
	\end{lemma}
	Then we come back to the last term of the left-hand side of (\ref{VFPre}). 
	Direct calculations show
	\begin{align*}
		T\left( \begin{array}{ccc}
			u^s|_{\Gamma_b}\\
			v^s|_{\Gamma_b}
		\end{array}\right) 
		&=\sum_{n\in\mathbb{Z}}M_n
		\left( \begin{array}{ccc}
			u_n\\
			v_n
		\end{array}\right) e^{i\alpha_n x_1+i\beta_n b}
		=\sum_{n\in\mathbb{Z}}\frac{1}{\kappa^2}
		\left(\begin{array}{ccc}   
			-i\omega\eps\beta_n u_n+ i\gamma \alpha_n v_n \\  
			-i\gamma \alpha_n u_n -i\omega\mu\beta_n v_n \\  
		\end{array}
		\right)e^{i\alpha_n x_1+i\beta_n b} ,\\
		T\left( \begin{array}{ccc}
			u^i|_{\Gamma_b}\\
			v^i|_{\Gamma_b}
		\end{array}\right) 
		&=M_0\left( \begin{array}{ccc}
			p_3\\
			q_3
		\end{array}\right)e^{i\alpha x_1-i\beta b}
		=-\frac{1}{\kappa^2}
		\left(\begin{array}{ccc}   
			i\omega\eps\beta p_3- i\gamma \alpha q_3 \\  
			i\gamma \alpha p_3+ i\omega\mu\beta q_3\\  
		\end{array}
		\right)e^{i\alpha x_1-i\beta b}.
	\end{align*}
	
	Therefore,
	\begin{align}
		&\quad \frac{1}{\kappa^2}
		\left(\begin{array}{ccc}   
			\omega\eps\partial_\nu u+ \gamma \partial_\tau v \\  
			\omega\mu\partial_\nu v- \gamma \partial_\tau u \\  
		\end{array}
		\right)\bigg|_{\Gamma_b}\nonumber\\
		&=\sum_{n\in\mathbb{Z}}	\frac{1}{\kappa^2}\left(                 
		\begin{array}{ccc}   
			i\omega\eps\beta_n & -i\gamma \alpha_n \\  
			i\gamma \alpha_n & i\omega\mu\beta_n \\  
		\end{array}
		\right) 
		\left( \begin{array}{ccc}
			u_n\\
			v_n
		\end{array}\right) e^{i\alpha_n x_1+i\beta_n b}
		-\frac{1}{\kappa^2}
		\left(\begin{array}{ccc}   
			i\omega\eps\beta p_3+ i\gamma \alpha q_3 \\  
			i\omega\mu\beta q_3-i\gamma \alpha p_3 \\  
		\end{array}
		\right)e^{i\alpha x_1-i\beta b}\nonumber\\
		&=-T\left( \begin{array}{ccc}
			u^s|_{\Gamma_b}\\
			v^s|_{\Gamma_b}
		\end{array}\right)
		-\frac{1}{\kappa^2}
		\left(\begin{array}{ccc}   
			i\omega\eps\beta p_3+ i\gamma \alpha q_3 \\  
			i\omega\mu\beta q_3-i\gamma \alpha p_3 \\  
		\end{array}
		\right)e^{i\alpha x_1-i\beta b}\nonumber\\
		&=-T\left( \begin{array}{ccc}
			u|_{\Gamma_b}\\
			v|_{\Gamma_b}
		\end{array}\right)
		+T\left( \begin{array}{ccc}
			u^i|_{\Gamma_b}\\
			v^i|_{\Gamma_b}
		\end{array}\right)
		-\frac{1}{\kappa^2}
		\left(\begin{array}{ccc}   
			i\omega\eps\beta p_3+ i\gamma \alpha q_3 \\  
			i\omega\mu\beta q_3-i\gamma \alpha p_3 \\  
		\end{array}
		\right)e^{i\alpha x_1-i\beta b}\nonumber\\
		&=-T\left( \begin{array}{ccc}
			u|_{\Gamma_b}\\
			v|_{\Gamma_b}
		\end{array}\right)
		-\frac{1}{\kappa^2}
		\left(\begin{array}{ccc}   
			i\omega\eps\beta p_3- i\gamma \alpha q_3 \\  
			i\gamma \alpha p_3 +i\omega\mu\beta q_3\\  
		\end{array}
		\right)e^{i\alpha x_1-i\beta b}
		-\frac{1}{\kappa^2}
		\left(\begin{array}{ccc}   
			i\omega\eps\beta p_3+ i\gamma \alpha q_3 \\  
			i\omega\mu\beta q_3-i\gamma \alpha p_3 \\  
		\end{array}
		\right)e^{i\alpha x_1-i\beta b}\nonumber\\
		&=-T\left( \begin{array}{ccc}
			u|_{\Gamma_b}\\
			v|_{\Gamma_b}
		\end{array}\right)
		-\frac{2}{\kappa^2}
		\left(\begin{array}{ccc}   
			i\omega\eps\beta p_3\\  
			i\omega\mu\beta q_3 \\  
		\end{array}
		\right)e^{i\alpha x_1-i\beta b}. \label{VFPreLsterm}
	\end{align}
	Note that in deriving (\ref{VFPreLsterm}), we have used the expression of $(u,v)$ given by 
	\begin{align*}
		u &= p_{3}e^{i\alpha x_{1} - i\beta x_{2}} + \sum_{n \in \mathbb{Z}} u_{n}e^{i\alpha_{n}x_{1} + i\beta_{n}x_{2}}, \\
		v &= q_{3}e^{i\alpha x_{1} - i\beta x_{2}} + \sum_{n \in \mathbb{Z}} v_{n}e^{i\alpha_{n}x_{1} + i\beta_{n}x_{2}},\qquad  x_2>\Gamma_{\max}.
	\end{align*}
	Inserting (\ref{VFPreLsterm}) into (\ref{VFPre}), we get the variational formulation 
	\begin{align}
		B(u,v;\varphi,\psi)=F(\varphi,\psi) \qquad \text{for all } \,(\varphi,\psi) \in X,\label{VF}
	\end{align}
	where
	\begin{align}
		B(u,v;\varphi,\psi)&:=
		\int_{\Gamma}\frac{i}{\lambda}u\,\overline{\varphi}+i\lambda v\,\overline{\psi}\,ds+
		\int_{\Omega_b}\left[\frac{\omega\eps}{\kappa^2}\nabla u\cdot \nabla\overline{\varphi}-\frac{\gamma}{\kappa^2}\nabla v\cdot \nabla^\bot\overline{\varphi}+\frac{\omega\mu}{\kappa^2}\nabla v\cdot \nabla\overline{\psi}\right.\nonumber\\
		&\quad+\left.\frac{\gamma}{\kappa^2}\nabla u\cdot \nabla^\bot\overline{\psi}-\omega \eps\,u\;\overline{\varphi}-\omega \mu\,v\;\overline{\psi}\right]\;dx 
		+\int_{\Gamma_b}
		T\left( \begin{array}{ccc}
			u\\
			v
		\end{array}\right)
		\cdot
		\left( \begin{array}{ccc}
			\overline{\varphi}\\
			\overline{\psi}
		\end{array}\right)\,ds,\label{VFB}\\
		F(\varphi,\psi)&:=-\frac{2i\omega\eps\beta e^{-i\beta b}}{\kappa^2}\int_{\Gamma_b}(\eps p_3\overline{\varphi}+\mu q_3\overline{\psi})e^{i\alpha x_1}\,ds.\label{VFF}
	\end{align}
	Below we prove an energy formula under the impedance boundary condition. 
	\begin{lemma}
		Let $u,v\in H_\alpha^1(\Omega_b)$ be the total fields to our conical diffraction problem. 
		We have the energy formula
		\begin{align}
			\frac{2\pi\omega}{\kappa^2}\sum_{|\alpha_n|\leq \kappa} \beta_n(\eps |u_n|^2+\mu |v_n|^2)=\int_{\Gamma}\frac{1}{\lambda}|u|^2+\lambda |v|^2\,ds+\frac{2\pi\omega\beta}{\kappa^2}(\eps|p_3|^2+\mu|q_3|^2). \label{eqn:energyF}
		\end{align}
		\begin{proof}
			By $ (\ref{VFPre}) $ and taking $ \varphi=u,\psi=v $, we have
			\begin{align}
				0&=\int_{\Gamma}\frac{i}{\lambda}|u|^2+i\lambda |v|^2\,ds-\frac{1}{\kappa^2}\int_{\Gamma_b}\begin{pmatrix}
					\omega\eps\partial_n u+ \gamma \partial_\tau v \\
					\omega\mu\partial_n v- \gamma \partial_\tau u
				\end{pmatrix}
				\cdot 
				\begin{pmatrix}
					\overline{u} \\
					\overline{v}
				\end{pmatrix} \,ds \nonumber\\
				&\quad+\int_{\Omega_b}\left[\frac{\omega\eps}{\kappa^2}|\nabla u|^2-\frac{\gamma}{\kappa^2}\nabla v\cdot \nabla^\bot\overline{u}-\omega \eps|u|^2+\frac{\omega\mu}{\kappa^2}|\nabla v|^2+\frac{\gamma}{\kappa^2}\nabla u\cdot \nabla^\bot\overline{v}-\omega \mu|v|^2\right]\;dx .\label{energy1}
			\end{align}
			We want to calculate the imaginary part of (\ref{energy1}). 
			First, we have
			\begin{align*}
				&\quad	\mbox{Im}\,\int_{\Gamma_b} \begin{pmatrix}
					\omega\eps\partial_n u+ \gamma \partial_\tau v \\
					\omega\mu\partial_n v- \gamma \partial_\tau u
				\end{pmatrix}
				\cdot 
				\begin{pmatrix}
					\overline{u} \\
					\overline{v}
				\end{pmatrix}\,ds\nonumber\\
				&=	\mbox{Im}\,\int_{0}^{2\pi}\left(\begin{array}{ccc}   
					\omega\eps\partial_2 u- \gamma \partial_1 v \\  
					\omega\mu\partial_2 v+ \gamma \partial_1 u \\  
				\end{array}
				\right)\bigg|_{\Gamma_b}\cdot
				\left( \begin{array}{ccc}
					\overline{u}\\
					\overline{v}
				\end{array}\right)\bigg|_{\Gamma_b}\,dx_1
			\end{align*}
			For the total field $ u $ and $ v $, we can easily get
			\begin{align*}
				&\quad	\mbox{Im}\,\int_{0}^{2\pi}\left[\omega\eps(-i\beta p_3e^{i\alpha x_1-i\beta b}+\sum_{n\in\mathbb{Z}}i\beta_n u_ne^{i\alpha_n x_1+i\beta_n b})-\gamma(i\alpha q_3e^{i\alpha x_1-i\beta b}\right.\nonumber\\
				&\quad\left.+\sum_{n\in\mathbb{Z}}i\alpha_n v_ne^{i\alpha_n x_1+i\beta_n b})\right]\left(\overline{p_3e^{i\alpha x_1-i\beta b}+\sum_{m\in\mathbb{Z}}u_m e^{i\alpha_m x_1+i\beta_m b}}\right)\nonumber\\
				&\quad+ \left[\omega \mu(-i\beta q_3e^{i\alpha x_1-i\beta b}+\sum_{n\in\mathbb{Z}}i\beta_n v_ne^{i\alpha_n x_1+i\beta_n b})+\gamma(i\alpha p_3e^{i\alpha x_1-i\beta b}\right.\nonumber\\
				&\quad \left.+\sum_{n\in\mathbb{Z}}i\alpha_n u_ne^{i\alpha_n x_1+i\beta_n b})\right]\left(\overline{q_3e^{i\alpha x_1-i\beta b}+\sum_{m\in\mathbb{Z}}v_m e^{i\alpha_m x_1+i\beta_m b}}\right)\,dx_1\nonumber\\
				&=\mbox{Im}\, 2\pi \left[\omega\eps\left(-i\beta |p_3|^2+\sum_{n\in\mathbb{Z}}i\beta_n |u_n|^2\right)-\gamma\left(i\alpha q_3\bar{p}_3+\sum_{n\in\mathbb{Z}}i\alpha_n v_n\bar{u}_n\right)\right.\nonumber\\
				&\quad+ \left.\omega \mu\left(-i\beta |q_3|^2+\sum_{n\in\mathbb{Z}}i\beta_n |v_n|^2\right)+\gamma\left(i\alpha p_3\bar{q}_3+\sum_{n\in\mathbb{Z}}i\alpha_n u_n\bar{v}_n\right)\right].
			\end{align*}
			Therefore,
			\begin{align}
				&\quad 	\mbox{Im}\,\int_{\Gamma_b} \begin{pmatrix}
					\omega\eps\partial_n u+ \gamma \partial_\tau v \\
					\omega\mu\partial_n v- \gamma \partial_\tau u
				\end{pmatrix}
				\cdot 
				\begin{pmatrix}
					\overline{u} \\
					\overline{v}
				\end{pmatrix}\,ds \nonumber\\
				&= -2\pi\omega\beta\left(\eps|p_3|^2+\mu|q_3|^2\right)+2\pi\sum_{|\alpha_n|\leq \kappa} \omega\beta_n\left(\eps |u_n|^2+\mu |v_n|^2\right).\label{energy2}
			\end{align}
			In addition,
			\begin{align}
				&\quad	\mbox{Im}\,\int_{\Gamma_b} \nabla v\cdot \nabla^\bot\overline{u}-\nabla u\cdot \nabla^\bot\overline{v}\;dx \nonumber\\
				&=\mbox{Im}\,\int_{\Gamma_b}-\partial_1 v\partial_2 \bar{u}+\partial_2 v\partial_1 \bar{u}-(-\partial_1 u\partial_2 \bar{v}+\partial_2 u\partial_1 \bar{v})\;dx \nonumber\\
				&=\mbox{Im}\,\int_{\Gamma_b} -(\partial_1 v\partial_2 \bar{u}+\partial_2 u\partial_1 \bar{v})+(\partial_2 v\partial_1 \bar{u}+\partial_1 u\partial_2 \bar{v})\;dx \nonumber\\
				&=0.\label{energy3}
			\end{align}
			Taking the imaginary part of $ (\ref{energy1}) $ and using $ (\ref{energy2}) $ and $ (\ref{energy3}) $, we obtain
			\begin{align*}
				0&=\int_{\Gamma}\frac{1}{\lambda}|u|^2+\lambda |v|^2\,ds 
				-\mbox{Im}\,\frac{1}{\kappa^2}\int_{0}^{2\pi}\left(\begin{array}{ccc}   
					\omega\eps\partial_2 u- \gamma \partial_1 v \\  
					\omega\mu\partial_2 v+ \gamma \partial_1 u \\  
				\end{array}
				\right)\bigg|_{\Gamma_b}\cdot
				\left( \begin{array}{ccc}
					\overline{u}\\
					\overline{v}
				\end{array}\right)\bigg|_{\Gamma_b}\,dx_1\nonumber\\
				&=\int_{\Gamma}\frac{1}{\lambda}|u|^2+\lambda |v|^2\,ds+\frac{2\pi\omega\beta}{\kappa^2}\left(\eps|p_3|^2+\mu|q_3|^2\right)-	\frac{2\pi\omega}{\kappa^2}\sum_{|\alpha_n|\leq \kappa} \beta_n\left(\eps |u_n|^2+\mu |v_n|^2\right),
			\end{align*}
			which completes the proof.
		\end{proof}
	\end{lemma}
	\begin{theorem}\label{TH1}
		Suppose that $\Gamma$ is a Lipschitz curve, $k^2\neq\gamma^2$ and the impedance coefficient $\lambda<0$. 
		Then, the variational problem (\ref{VF}) has at most one solution $(u,v)\in X$.
	\end{theorem}
	\begin{proof}
		To prove uniqueness, we assume $ u^i=v^i=0 $, i.e. $ p_3=q_3=0 $. 
		Choosing $ \varphi=u,\psi=v $ in  $ (\ref{VF}) $ and taking the imaginary part, we have
		\begin{align}
			\int_{\Gamma}\frac{1}{\lambda}|u|^2+\lambda |v|^2\,ds+\mbox{Im}\,\int_{\Gamma_b}
			T\left( \begin{array}{ccc}
				u\\
				v
			\end{array}\right)
			\cdot
			\left( \begin{array}{ccc}
				\overline{u}\\
				\overline{v}
			\end{array}\right)\,ds=0.\label{ImVF}
		\end{align}
		Next, we calculate the second term of $ (\ref{ImVF}) $. 
		By the definition of $ T $ (see Definition \ref{dtn}), we have for $w = (u, v)^\top$ that
		\begin{align*}
			\mbox{Im}\,\int_{\Gamma_b}Tw\cdot \overline{w}\,ds&=\mbox{Im}\,\int_{\Gamma_b}\sum_{n\in\mathbb{Z}}M_n\hat{w}_n e^{i\alpha_n x}\cdot\overline{\sum_{m\in\mathbb{Z}}\hat{w}_m e^{i\alpha_m x}}\,ds\\
			&=\mbox{Im}\, 2\pi\sum_{n\in\mathbb{Z}}M_n\hat{w}_n\cdot\overline{\hat{w}_n}\\
			&=2\pi\sum_{n\in\mathbb{Z}}(\mbox{Im}\,M_n)\hat{w}_n\cdot\overline{\hat{w}_n},
		\end{align*}
		where $\hat{w}_{n} = (\hat{w}_{n, 1}, \hat{w}_{n, 2}) = (u_{n}, v_{n})$.
		Recalling the expression of $ M_n $, we have
		\begin{align*}
			\mbox{Im}\,M_n&=\frac{1}{2i}(M_n-M_n^*)\\
			&=\frac{1}{2i}\frac{1}{\kappa^2}\left[ 
			\left(                 
			\begin{array}{ccc}   
				-i\omega\eps\beta_n & i\gamma \alpha_n \\  
				-i\gamma \alpha_n & -i\omega\mu\beta_n \\  
			\end{array}
			\right) 
			-\left(                 
			\begin{array}{ccc}   
				i\omega\eps\overline{\beta_n} & i\gamma \alpha_n \\  
				-i\gamma \alpha_n & i\omega\mu\overline{\beta_n} \\  
			\end{array}
			\right) 
			\right] \\
			&=\frac{1}{\kappa^2}
			\left(                 
			\begin{array}{ccc}   
				\mbox{Im}\,(-i\omega\eps\beta_n) & 0\\  
				0 & \mbox{Im}\,(-i\omega\mu\beta_n) \\  
			\end{array}
			\right) \\
			&=\begin{cases}
				\frac{1}{\kappa^2}
				\left(                 
				\begin{array}{ccc}   
					-\omega\eps\beta_n & 0\\  
					0 & -\omega\mu\beta_n \\  
				\end{array}
				\right)
				,&|\alpha_n|\leq \kappa,\\
				0,&|\alpha_n|> \kappa.
			\end{cases}
		\end{align*}
		Therefore,
		\begin{align*}
			\mbox{Im}\,\int_{\Gamma_b}Tw\cdot \overline{w}\,ds
			&=2\pi\sum_{|\alpha_n|\leq \kappa} \frac{1}{\kappa^2}
			\left(                 
			\begin{array}{ccc}   
				-\omega\eps\beta_n & 0\\  
				0 & -\omega\mu\beta_n \\  
			\end{array}
			\right)
			\left( \begin{array}{ccc}
				\hat{w}_{n1}\\
				\hat{w}_{n2}
			\end{array}\right)\cdot
			\left( \begin{array}{ccc}
				\overline{\hat{w}_{n1}}\\
				\overline{\hat{w}_{n2}}
			\end{array}\right)\\
			&=-\frac{2\pi\omega}{\kappa^2}\sum_{|\alpha_n|\leq \kappa} \beta_n(\eps |\hat{w}_{n1}|^2+\mu |\hat{w}_{n2}|^2)\leq 0.
		\end{align*}
		Inserting these results into $ (\ref{ImVF}) $, we have
		\begin{align*}
			\int_{\Gamma}\frac{1}{\lambda}|u|^2+\lambda |v|^2\,ds-\frac{2\pi\omega}{\kappa^2}\sum_{|\alpha_n|\leq \kappa} \beta_n(\eps |u_n|^2+\mu |v_n|^2)=0. 
		\end{align*}
		Noting that $ \lambda<0 $, we have $ u=v=0\;\mbox{on}\;\Gamma $. 
		By the impedance radiation condition (\ref{eqn:BC}), we have $ \partial_n u=\partial_n v=0\;\mbox{on}\;\Gamma$. 
		By Holmgren theorem,  $ u=v=0\;\mbox{in}\;\Omega$.
	\end{proof}
	\begin{remark}\label{RmkUnique}
		We can also prove the uniqueness result by taking the imaginary part of the energy formula (\ref{eqn:energyF}) with $ p_3=q_3=0 $.
	\end{remark}
	The proof of Theorem \ref{TH1} provides an alternative approach to the proof of the energy formula via matrix operations. 
	\begin{definition}[Strong ellipticity] \label{DefstronglyEllip}
		We call a bounded sesquilinear form $ B(\cdot,\cdot) $ given on some Hilbert space $ X $  strongly elliptic if there exists a complex number $ \theta $, $ |\theta|=1 $ and a compact form $ q(\cdot,\cdot) $ such that 
		\ben
	{	\rm{Re}}\,(\theta B(u,u))\geq c\|u\|_X^2-q(u,u) \qquad \text{ for all } u \in X,
		\enn
		for some constant $c > 0$.
	\end{definition}
	The following theorem establishes the strong ellipticity of the form $ (\ref{VFB}) $ and leads, together with Theorem $ \ref{TH1} $ or Remark $ \ref{RmkUnique} $, to the solvability results for the conical diffraction problem.
	\begin{theorem}\label{StronglyEllip}
		The sesquilinear form $ B $ defined in $ (\ref{VF}) $ is strongly elliptic over $ X $.
	\end{theorem}
	We divide the proof of Theorem $ \ref{StronglyEllip} $ into several lemmas. 
	It is convenient to reformulate the variational form $ (\ref{VFB}) $ as follows (see \cite{El00})
	\ben
	B(u,v;\varphi,\psi)=A(u,v;\varphi,\psi)+B_1(u,v;\varphi,\psi)+C(u,v;\varphi,\psi)+D(u,v;\varphi,\psi),
	\enn
	where
	\ben
	A(u,v;\varphi,\psi)&:=&\int_{\Gamma}\frac{i}{\lambda}u\,\overline{\varphi}+i\lambda v\,\overline{\psi}\,ds,\\
	C(u,v;\varphi,\psi)&:=&\int_{\Omega_b}\omega \eps\,u\;\overline{\varphi}+\omega \mu\,v\;\overline{\psi}\;dx,\\
	D(u,v;\varphi,\psi)&:=&\int_{\Gamma_b}
	T\left( \begin{array}{ccc}
		u\\
		v
	\end{array}\right)
	\cdot
	\left( \begin{array}{ccc}
		\overline{\varphi}\\
		\overline{\psi}
	\end{array}\right)\,ds,
	\enn
	and
	\ben
	B_1(u,v;\varphi,\psi)&:=&
	\int_{\Omega_b}\left[\frac{\omega\eps}{\kappa^2}\nabla u\cdot \nabla\overline{\varphi}-\frac{\gamma}{\kappa^2}\nabla v\cdot \nabla^\bot\overline{\varphi}+\frac{\omega\mu}{\kappa^2}\nabla v\cdot \nabla\overline{\psi}+\frac{\gamma}{\kappa^2}\nabla u\cdot \nabla^\bot\overline{\psi}\right]\;dx \nonumber\\
	&=&\int_{\Omega_b} \mathcal{D}(\partial_1 u,\partial_1 v,\partial_2 u,\partial_2 v)^\mathsf{T}\cdot \overline{(\partial_1 u,\partial_1 v,\partial_2 u,\partial_2 v)^\mathsf{T}}\;dx,
	\enn
	with the matrix $ \mathcal{D} $ given by (see \cite{El00})
	$$
	\mathcal{D}=\frac{1}{\kappa^2}
	\left(
	\begin{matrix}
		\omega \epsilon & 0 & 0 & -\gamma      \\
		0 & \omega\mu & \gamma & 0      \\
		0 & \gamma & \omega \epsilon & 0 \\
		-\gamma & 0 & 0 & \omega\mu      \\
	\end{matrix}
	\right).
	$$
	
	We can further write $B_1 $ into the matrix form
	\be \label{B1}
	B_1(u,v;\varphi,\psi)=\int_{\Omega_b}
	N^+ \partial^+\left( \begin{array}{ccc}
		u\\
		v
	\end{array}\right)
	\cdot \overline{\partial^+
		\left( \begin{array}{ccc}
			\varphi\\
			\psi
		\end{array}\right)}
	+ N^- \partial^-\left( \begin{array}{ccc}
		u\\
		v
	\end{array}\right)
	\cdot \overline{\partial^-
		\left( \begin{array}{ccc}
			\varphi\\
			\psi
		\end{array}\right)}\,ds
	\en
	where 
	\ben 
	N^{\pm}=\frac{1}{\kappa^2}\left(
	\begin{matrix}
		\omega \epsilon & \pm i\gamma \\
		\mp i\gamma & \omega\mu
	\end{matrix}
	\right),\quad
	\partial^+ :=\frac{1}{\sqrt{2}} (-i\partial_1+\partial_2),\quad \partial^- :=\frac{1}{\sqrt{2}} (\partial_1-i\partial_2).
	\enn
	To study the form $ B $, we need the following lemma.
	\begin{lemma}\label{PD}
		Choose $ \theta=\frac{i+\delta}{|i+\delta|} $ with $ \delta >0$ sufficiently small.\\
		\text{(i)} For any $ \xi\in \mathbb{C}^2 $, we have $ \mbox{Re}\,(\theta N^\pm \xi\cdot \bar{\xi})\geq C_N|\xi|^2 $, where 
		\be 
		C_N=\frac{1}{2\omega\epsilon\mu\cos^2\phi}{\rm{Re}}\,\theta \left[(\epsilon+\mu)-\sqrt{(\epsilon-\mu)^2+4\epsilon\mu\sin^2\phi}\right] \geqslant 0 .\label{cn}
		\en
		\text{(ii)} Let $M_n \in \mathbb{C}^{2\times 2}$ be defined by (\ref{Mn}). It holds that $\mbox{Re}\,(\theta M_n)\geq 0 $ for all $n \in \mathbb{Z}\backslash \mathcal{A}$, where the index set $\mathcal{A}$ is defined by
		\begin{equation}
			\begin{aligned}
				\mathcal{A} = \{n\in \mathbb{Z}:-k(1+\sin\theta\cos\phi)< n \leq -k\cos\phi(1+\sin\theta) \\
				\text{ or }\; k \cos\phi (1-\sin\theta) \leq n < k(1-\sin\theta\cos\phi)\}.
			\end{aligned}
			\label{aA}
		\end{equation}
	\end{lemma}
	\begin{proof}
		\text{(i)} By the definition of $ N^\pm $, we have
		\ben
		\mbox{Re}\,(\theta N^\pm)=\frac{\theta N^\pm+(\theta N^\pm)^*}{2}=\frac{1}{\kappa^2}
		\left(
		\begin{matrix}
			\omega \epsilon \mbox{Re}\,\theta & \pm i\gamma  \mbox{Re}\,\theta    \\
			\mp i\gamma  \mbox{Re}\,\theta & \omega\mu \mbox{Re}\,\theta      \\
		\end{matrix}
		\right),
		\enn
		which is a Hermitian matrix.  
		Recalling that $ \gamma=k\sin \phi=\omega\sqrt{\epsilon\mu}\sin \phi $ and $ \kappa^2=k^2 \cos^2\phi$, we compute the eigenvalues of $ \mbox{Re}\,(\theta N^\pm) $ as following
		\ben
		\lambda_{1}&=&\frac{1}{2\omega\epsilon\mu\cos^2\phi}\mbox{Re}\,\theta\left[(\epsilon+\mu)+\sqrt{(\epsilon-\mu)^2+4\epsilon\mu\sin^2\phi}\right] > 0,\\
		\lambda_{2}&=&\frac{1}{2\omega\epsilon\mu\cos^2\phi}\mbox{Re}\,\theta\left[[(\epsilon+\mu)-\sqrt{(\epsilon-\mu)^2+4\epsilon\mu\sin^2\phi}\right] \geq 0.
		\enn
		Defining $ C_N=\lambda_2 < \lambda_1$, and by \cite[Theorem 4.2.2]{HornMatrixAnalysis}, we complete the first part of the proof.\\
		\text{(ii)} Recalling the definition of $ M_n $, we have
		\ben
		\theta M_n= \frac{1}{\kappa^2}
		\left(                 
		\begin{array}{ccc}   
			-i\omega\eps\beta_n \theta & i\gamma \alpha_n\theta\\  
			-i\gamma \alpha_n\theta & -i\omega\mu\beta_n \theta \\  
		\end{array}
		\right)=
		\frac{1}{\kappa^2}\left(                 
		\begin{array}{ccc}   
			-i(\omega \mu)^{-1}k^2\beta_n\theta  & i\gamma \alpha_n\theta\\  
			-i\gamma \alpha_n\theta &-i\omega \mu \beta_n \theta \\  
		\end{array}
		\right).
		\enn
		\text{Case 1.} $ |\alpha_n|< \kappa $, i.e., $ \beta_n \in \mathbb{R}$ is real number. 
		We have 
		$$ (\theta M_n)^*=\frac{1}{\kappa^2} \left(                 
		\begin{array}{ccc}   
			i(\omega \mu)^{-1}k^2\beta_n \bar{\theta}  & i\gamma \alpha_n\bar{\theta}\\  
			-i\gamma \alpha_n\bar{\theta} & i\omega \mu \beta_n \bar{\theta} \\  
		\end{array}
		\right). $$
		Therefore,
		\ben
		\mbox{Re}\,(\theta M_n)=\frac{\theta M_n+(\theta M_n)^*}{2}=\frac{1}{\kappa^2}\left(                 
		\begin{array}{ccc}   
			(\omega \mu)^{-1}k^2\beta_n \mbox{Im}\,\theta  & i\gamma \alpha_n\mbox{Re}\,\theta\\  
			-i\gamma \alpha_n\mbox{Re}\,\theta & \omega \mu \beta_n \mbox{Im}\,\theta \\  
		\end{array}
		\right).
		\enn
		In this case, $ \mbox{Re}\,(\theta M_n)\geq 0 $ if and only if the following two conditions are satisfied:
		\be
		\mbox{Im}\,\theta &\geq& 0,\nonumber\\
		det(\mbox{Re}\,(\theta M_n))&=&\frac{1}{\kappa^4}\left[(\omega \mu)^{-1}k^2\beta_n \mbox{Im}\,\theta\omega \mu \beta_n \mbox{Im}\,\theta-\gamma^2\alpha_n^2(\mbox{Re}\,\theta)^2\right]\nonumber\\
		&=&\frac{1}{\kappa^4}\left[k^2\beta_n^2-\gamma^2\alpha_n^2\delta^2\right](\mbox{Im}\,\theta)^2\geq 0.\label{cs1}
		\en
		The conditions in (\ref{cs1}) obviously hold due to the definition of $ \theta $ with a small $  \delta>0$.
		
		\text{Case 2.} $ |\alpha_n|\geq \kappa $, i.e., $ \beta_n $ is a pure imaginary number. 
		We have 
		$$ (\theta M_n)^*= \frac{1}{\kappa^2}\left(                 
		\begin{array}{ccc}   
			(\omega \mu)^{-1}k^2|\beta_n| \bar{\theta}  & i\gamma \alpha_n\bar{\theta}\\  
			-i\gamma \alpha_n\bar{\theta} & \omega \mu |\beta_n| \bar{\theta} \\  
		\end{array}
		\right). $$
		Therefore,
		\ben
		\mbox{Re}\,(\theta M_n)=\frac{\theta M_n+(\theta M_n)^*}{2}=\frac{1}{\kappa^2}\left(                 
		\begin{array}{ccc}   
			(\omega \mu)^{-1}k^2|\beta_n| \mbox{Re}\,\theta  & i\gamma \alpha_n\mbox{Re}\,\theta\\  
			-i\gamma \alpha_n\mbox{Re}\,\theta &\omega \mu |\beta_n| \mbox{Re}\,\theta \\  
		\end{array}
		\right).
		\enn
		In this case, $ \mbox{Re}\,(\theta M_n)\geq 0 $ if and only if the following two conditions are satisfied:
		\ben
		\mbox{Re}\,\theta &\geq& 0,\\
		det(\mbox{Re}\,(\theta M_n))&=&\frac{1}{\kappa^4}\left[(\omega \mu)^{-1}k^2|\beta_n| (\mbox{Re}\,\theta)\omega \mu |\beta_n| \mbox{Re}\,\theta-\gamma^2\alpha_n^2 (\mbox{Re}\,\theta)^2\right]\\
		&=&\frac{1}{\kappa^4}\left[k^2|\beta_n|^2-\gamma^2\alpha_n^2\right](\mbox{Re}\,\theta)^2\geq 0.
		\enn
		The first condition is obvious. The second condition can be fulfilled if $ k^2|\beta_n|^2-\gamma^2\alpha_n^2\geq 0 $. Recalling that $ |\beta_n|^2= \alpha_n^2-(k^2-\gamma^2)$, we have $ \alpha_n^2\geq k^2 $. 
		Combining the above two cases, we get that when 
		\ben
		n\in \mathcal{B}:=\{n\in\mathbb{Z}:k^2-\gamma^2\leq \alpha_n^2< k^2 \}, 
		\enn
		$ \mbox{Re}\,(\theta M_n) $ is not positive definite. 
		We should point out that $ n\in \mathcal{B}$ if $ \beta_ n= 0 $ (i.e. $|\alpha_n|=\kappa$).
		Next, we continue to simplify the set $ \mathcal{B} $. 
		Recalling $ \alpha=k\sin\theta\cos\phi $ and $ \alpha_n=n+ \alpha$, we have
		\begin{align*}
			\mathcal{B}=&\{n\in\mathbb{Z}:k^2-\gamma^2\leq \alpha_n^2< k^2 \}\\
			=&\{n\in\mathbb{Z}:k^2-k^2\sin^2\phi \leq (k \sin\theta\cos\phi+n)^2<k^2\}\\
			=&\{n\in\mathbb{Z}:k\cos\phi\leq |k \sin\theta\cos\phi+n|<k\}\\
			=&\{n\in\mathbb{Z}: -k(1+\sin\theta\cos\phi)<n< k(1-\sin\theta\cos\phi)\} \cap \\
			&\{n\in\mathbb{Z}: n \geq k\cos\phi(1-\sin\theta) \text{ or } n\leq -k\cos\phi(1+\sin\theta)\} \\
			=&\{n\in \mathbb{Z}:-k(1+\sin\theta\cos\phi)< n \leq -k\cos\phi(1+\sin\theta)\; \\
			& \text{ or } k \cos\phi (1-\sin\theta) \leq n < k(1-\sin\theta\cos\phi)\}.
		\end{align*}
		This coincides with the set $ \mathcal{A} $ given by (\ref{aA}).
	\end{proof}
	
	For $ u|_{\Gamma_b}=\sum_{n\in\mathbb{Z}}\tilde{u}_n e^{i\alpha_n x_1} ,\,v|_{\Gamma_b}=\sum_{n\in\mathbb{Z}}\tilde{v}_n e^{i\alpha_n x_1}$, we define 
	\begin{align}
		q(u,v;u,v)=
		2\pi\mbox{Re}\,\sum_{n\in \mathcal{A}}\theta M_n\left( \begin{array}{ccc}
			\tilde{u}_n\\
			\tilde{v}_n
		\end{array}\right)
		\cdot
		\left( \begin{array}{ccc}
			\overline{\tilde{u}_n}\\
			\overline{\tilde{v}_n}
		\end{array}\right),\label{defQ}
	\end{align}
	where the set $\mathcal{A}$ is defined by (\ref{aA}).
	
	\textbf{Proof of Theorem \ref{StronglyEllip}.}
	Choose $ \theta=\frac{i+\delta}{|i+\delta|} $. 
	By the definition of $\mathcal{A}$,
	\begin{align*}
		\mbox{Re}\,(\theta A(u,v;u,v))&=\mbox{Re}\,\int_{\Gamma}\frac{i+\delta}{|i+\delta|}\left(\frac{i}{\lambda}|u|^2+i\lambda |v|^2\right)\;ds\\
		&=-\int_{\Gamma}\frac{1}{|i+\delta|}\left(\frac{1}{\lambda}|u|^2+\lambda |v|^2\right)\;ds \geq 0.
	\end{align*}
	Before calculating $ \mbox{Re}\,(\theta B_1(u,v;u,v)) $, we compute the following relation: 
	\begin{align*}
		&\quad \partial^+\left( \begin{array}{ccc}
			u\\
			v
		\end{array}\right)
		\cdot \overline{\partial^+
			\left( \begin{array}{ccc}
				u\\
				v
			\end{array}\right)}=\frac{1}{2}
		\left( \begin{array}{ccc}
			-i\partial_1 u+\partial_2 u\\
			-i\partial_1 v+\partial_2 v
		\end{array}\right)
		\cdot
		\left( \begin{array}{ccc}
			i\partial_1 \bar{u}+\partial_2 \bar{u}\\
			i\partial_1 \bar{v}+\partial_2 \bar{v}
		\end{array}\right)\\
		&=\frac{1}{2}\left[|\partial_1 u|^2+|\partial_2 u|^2+i\partial_1 \bar{u}\partial_2 u-i\partial_1 u\partial_2 \bar{u}+|\partial_1 v|^2+|\partial_2 v|^2+i\partial_1 \bar{v}\partial_2 v-i\partial_1 v\partial_2 \bar{v}\right]\\
		&=\frac{1}{2}\left[|\nabla u|^2+|\nabla v|^2+i(\partial_1 \bar{u}\partial_2 u+\partial_1 \bar{v}\partial_2 v)-i(\partial_1 u\partial_2 \bar{u}+\partial_1 v\partial_2 \bar{v})\right].
	\end{align*}
	Similarly,
	\begin{align*}
		\partial^-\left( \begin{array}{ccc}
			u\\
			v
		\end{array}\right)
		\cdot \overline{\partial^-
			\left( \begin{array}{ccc}
				u\\
				v
			\end{array}\right)}
		=\frac{1}{2}\left[|\nabla u|^2+|\nabla v|^2+i(\partial_1 u\partial_2 \bar{u}+\partial_1 v\partial_2 \bar{v})-i(\partial_1 \bar{u}\partial_2 u+\partial_1 \bar{v}\partial_2 v)\right].
	\end{align*}
	Therefore,
	\begin{align*}
		\partial^+\left( \begin{array}{ccc}
			u\\
			v
		\end{array}\right)
		\cdot \overline{\partial^+
			\left( \begin{array}{ccc}
				u\\
				v
			\end{array}\right)}
		+  \partial^-\left( \begin{array}{ccc}
			u\\
			v
		\end{array}\right)
		\cdot \overline{\partial^-
			\left( \begin{array}{ccc}
				u\\
				v
			\end{array}\right)}=|\nabla u|^2+|\nabla v|^2.
	\end{align*}
	Then, by (\ref{B1}) and Lemma \ref{PD}, we have
	\begin{align*}
		&\quad\mbox{Re}\,(\theta B_1(u,v;u,v))\\
		&=\mbox{Re}\,\left(\theta \int_{\Omega_b}
		N^+ \partial^+\left( \begin{array}{ccc}
			u\\
			v
		\end{array}\right)
		\cdot \overline{\partial^+
			\left( \begin{array}{ccc}
				u\\
				v
			\end{array}\right)}
		+ N^- \partial^-\left( \begin{array}{ccc}
			u\\
			v
		\end{array}\right)
		\cdot \overline{\partial^-
			\left( \begin{array}{ccc}
				u\\
				v
			\end{array}\right)}\,dx\right)\\
		&\geq C_N \int_{\Omega_b}
		\partial^+\left( \begin{array}{ccc}
			u\\
			v
		\end{array}\right)
		\cdot \overline{\partial^+
			\left( \begin{array}{ccc}
				u\\
				v
			\end{array}\right)}
		+  \partial^-\left( \begin{array}{ccc}
			u\\
			v
		\end{array}\right)
		\cdot \overline{\partial^-
			\left( \begin{array}{ccc}
				u\\
				v
			\end{array}\right)}\,dx \\
		&=C_N \int_{\Omega_b}|\nabla u|^2+|\nabla v|^2 \,dx\\
		&=C_N\left(\|u\|^2_{H^1(\Omega_b)}+\|v\|^2_{H^1(\Omega_b)}\right)-C_N \int_{\Omega_b} | u|^2+| v|^2 \,dx,	
	\end{align*}
	where $C_{N} \geq 0$ is defined by (\ref{cn}).
	It is obvious that
	\begin{align*}
		\mbox{Re}\,(\theta C(u,v;u,v))&=\mbox{Re}\,\int_{\Omega_b}\frac{i+\delta}{|i+\delta|}\left(\omega \eps|u|^2+\omega \mu|v|^2\right)\;dx\\
		&=\frac{\delta}{|i+\delta|}\int_{\Omega_b}\omega \eps|u|^2+\omega \mu|v|^2\;dx \geq 0.
	\end{align*}
	Then we need to consider the term $ B(u,v;\varphi,\psi) $. 
	Suppose that
	\begin{equation*}
		u|_{\Gamma_b}=\sum_{n\in\mathbb{Z}}\tilde{u}_n e^{i\alpha_n x_1} ,\quad\,v|_{\Gamma_b}=\sum_{n\in\mathbb{Z}}\tilde{v}_n e^{i\alpha_n x_1}.
	\end{equation*}
	We get
	\begin{align*}
		\mbox{Re}\,(\theta D(u,v;u,v))&=
		\mbox{Re}\,\left(\theta\int_{\Gamma_b}
		T\left( \begin{array}{ccc}
			u\\
			v
		\end{array}\right)
		\cdot
		\left( \begin{array}{ccc}
			\overline{u}\\
			\overline{v}
		\end{array}\right)\,ds\right)\\
		&=2\pi\mbox{Re}\,\left(\sum_{n\in\mathbb{Z}}\theta M_n 
		\left( \begin{array}{ccc}
			\tilde{u}_n\\
			\tilde{v}_n
		\end{array}\right)
		\cdot
		\left( \begin{array}{ccc}
			\overline{\tilde{u}_n}\\
			\overline{\tilde{v}_n}
		\end{array}\right)
		\right)\\
		&=2\pi\mbox{Re}\,\left(\sum_{n\in\mathbb{Z}/\mathcal{A}}\theta M_n 
		\left( \begin{array}{ccc}
			\tilde{u}_n\\
			\tilde{v}_n
		\end{array}\right)
		\cdot
		\left( \begin{array}{ccc}
			\overline{\tilde{u}_n}\\
			\overline{\tilde{v}_n}
		\end{array}\right)
		\right)+q(u,v;u,v) \\
		&\geq q(u,v;u,v),
	\end{align*}
	where we have used Lemma \ref{PD} (ii). 
	Note that $q(u,v;u,v)$ is a compact form because $\mathcal{A}$ is a finite set. 
	Therefore, by Lemma \ref{PD}, we have 
	\begin{align*}
		&\quad \mbox{Re}\,(\theta B(u,v;u,v))\\
		&=\mbox{Re}\,(\theta A(u,v;u,v))+\mbox{Re}\,(\theta B_1(u,v;u,v))+\mbox{Re}\,(\theta C(u,v;u,v))+\mbox{Re}\,(\theta D(u,v;u,v))\\
		&\geq C_N\left(\|u\|^2_{H^1(\Omega_b)}+\|v\|^2_{H^1(\Omega_b)}\right)-Q(u,v;u,v),
	\end{align*}
	where $C_{N} \geq 0$ is defined by (\ref{cn}) and
	\begin{align*}
		Q(u,v;u,v) := C_N\int_{\Omega_b}|u|^2+|v|^2\,dx-\frac{\delta}{|i+\delta|}\int_{\Omega_b}\omega \eps|u|^2+\omega \mu|v|^2\;dx-q(u,v;u,v)
	\end{align*}
	is a compact form over $X \times X$. By Definition \ref{DefstronglyEllip}, we finish the proof.
	$\hfill\Box$
	\begin{theorem}
		Suppose that $\Gamma$ is a Lipschitz curve, $k^2\neq\gamma^2$ and that the impedance coefficient $\lambda<0$. Then, the variational problem (\ref{v1})-(\ref{v2}) admits a unique solution $(u,v)\in X$.
	\end{theorem}
	\begin{proof}
		Under the assumption of Theorem \ref{StronglyEllip}, the operator defined in (\ref{VFB}) is a Fredholm operator with index zero. 
		Using Theorem \ref{TH1}, we obtain the existence and uniqueness result as a consequence of the Fredholm alternative.
	\end{proof}
	
	\section{Finite element analysis} \label{sec:FEM}
	
	We study the finite element approximation of the variational problem $ (\ref{VF}) $. 
	Let $ \{X_h^2: h \in (0,1)\} $ be a family of finite dimensional subspaces of $ H^1_{\alpha}(\Omega_b)^2 $, where $ h $ stands for the maximum mesh size after partitioning $ \Omega_b $ into simple domains, for example, a regular triangulation of $ \Omega_b $. 
	We make a general assumption \cite{Ciarlet2002} on the  subspace $ X_h^2 $ for 
	$ (\varphi,\psi) \in H^\rho_{\alpha}(\Omega_b)^2,\; \rho\geq 2 $,
	\begin{align}
		&\inf_{(\xi,\eta)\in X_h^2 } \left( \|(\varphi,\psi)-(\xi,\eta)\|_{L^2(\Omega_b)^2}+h\|(\nabla\varphi,\nabla\psi)-(\nabla\xi,\nabla\eta)\|_{L^2(\Omega_b)^2}\right.\nonumber\\
		&\qquad\qquad +h^{1/2}\|(\varphi,\psi)-(\xi,\eta)\|_{L^2(\Gamma_b)^2}
		+h\|(\varphi,\psi)-(\xi,\eta)\|_ {H^{1/2}(\Gamma_b)^2} \nonumber\\ &\qquad\qquad\left.+h^{1/2}\|(\varphi,\psi)-(\xi,\eta)\|_{L^2(\Gamma_b)^2} +h\|(\varphi,\psi)-(\xi,\eta)\|_ {H^{1/2}(\Gamma)^2}\right)\nonumber\\
		&\leq C h^l \|(\varphi,\psi)\|_{H^l(\Omega)^2},\quad l\in [2,\rho] \label{asuptnXH}
	\end{align}
	where the positive constant $ C $ is independent of $ h $ and $ (\varphi,\psi) $. 
	The finite element approximation to the variational $ (\ref{VF}) $ is to find $ (u_h,v_h) \in X_h^2 $ such that
	\be
	B(u_h,v_h;\varphi_h,\psi_h)=F(\varphi_h,\psi_h),\quad \text{ for all}  \;(\varphi_h,\psi_h)\in X_h^2, \label{VFFEM}
	\en
	where $ B $ is defined by $ (\ref{VFB}) $ and $ F $ is defined by $ (\ref{VFF}) $.
	The finite element method consists of the following steps to solve $ (\ref{VFFEM}) $: 
	\begin{description}
		\item[(1)] Choose a finite set of basis functions $ \{\phi_1,\phi_2,\cdots,\phi_m\} $ of $ X_h $;
		\item[(2)]Let $u_h=c_1 \phi_1+ c_2 \phi_2+\cdots+c_m \phi_m$, $v_h=d_1 \phi_1+ d_2 \phi_2+\cdots+d_m \phi_m$. Substitute the expression into $ (\ref{VFFEM}) $ and choose $ (\varphi_h,\psi_h)=(\phi_i,0),\; (0,\phi_i),\; i=1,2,\cdots,m$ to get a system of linear equations;
		\item[(3)] Solve the linear system for the coefficients $ c_1,c_2,\cdots,c_m,d_1,d_2,\cdots,d_m $ and get the approximation of $ (u,v) $ in $ X_h^2 $.
	\end{description}
	More precisely, we have
	\begin{align*}
		&\quad	B(u_h,v_h;\phi_i,0)\\
		&=
		\int_{\Omega_b}\left[\frac{\omega\eps}{\kappa^2}\left(\sum_{j=1}^{m}c_j \nabla\phi_j\right)\cdot \nabla\overline{\phi}_i-\frac{\gamma}{\kappa^2}\left(\sum_{j=1}^{m}d_j \nabla\phi_j\right)\cdot \nabla^\bot\overline{\phi}_i-\omega \eps\,\left(\sum_{j=1}^{m}c_j \phi_j\right)\;\overline{\phi}_i\right]\;dx \nonumber\\
		&\quad+\int_{\Gamma}\frac{i}{\lambda}\left(\sum_{j=1}^{m}c_j \phi_j\right)\,\overline{\phi}_i\,ds+\int_{\Gamma_b}
		T\left( \begin{array}{ccc}
			\sum_{j=1}^{m}c_j \phi_j\\
			\sum_{j=1}^{m}d_j \phi_j
		\end{array}\right)
		\cdot
		\left( \begin{array}{ccc}
			\overline{\phi}_i\\
			0
		\end{array}\right)\,ds,\\
	\end{align*}
	\begin{align*}
		&\quad	B(u_h,v_h;0,\phi_i)\\
		&=
		\int_{\Omega_b}\left[\frac{\omega\mu}{\kappa^2}\left(\sum_{j=1}^{m}d_j \nabla\phi_j\right)\cdot \nabla\overline{\phi}_i+\frac{\gamma}{\kappa^2}\left(\sum_{j=1}^{m}c_j \nabla\phi_j\right) \cdot \nabla^\bot\overline{\phi}_i-\omega \mu\,\left(\sum_{j=1}^{m}d_j \phi_j\right)\;\overline{\phi}_i\right]\;dx \nonumber\\
		&\quad+\int_{\Gamma}i\lambda \left(\sum_{j=1}^{m}d_j \phi_j\right)\;\overline{\phi}_i\,ds
		+\int_{\Gamma_b}
		T\left( \begin{array}{ccc}
			\sum_{j=1}^{m}c_j \phi_j\\
			\sum_{j=1}^{m}d_j \phi_j
		\end{array}\right)
		\cdot
		\left( \begin{array}{ccc}
			0\\
			\overline{\phi}_i
		\end{array}\right)\,ds.
	\end{align*}
	In order to deduce the stiffness matrix, we need to define the following inner product.
	\begin{align*}
		\langle f,g \rangle _{\Omega_b}=\int_{\Omega_b}f\bar{g}\;dx, \quad
		\langle f,g \rangle _{\Gamma}=\int_{\Gamma}f\bar{g}\;ds, \quad
		\langle f,g \rangle _{\Gamma_b}=\int_{\Gamma_b}f\bar{g}\;ds.
	\end{align*}
	Let $ \mathcal{B}\in \mathbb{C}^{2m \times 2m} $ be the stiffness matrix with the entries
	\begin{align*}
		B_{ij}= \begin{cases}
			\frac{\omega\eps}{\kappa^2}\langle \nabla \phi_j,\nabla \phi_i\rangle_{\Omega_b}-\omega\epsilon \langle  \phi_j, \phi_i\rangle_{\Omega_b}+\frac{i}{\lambda}\langle  \phi_j, \phi_i\rangle_{\Gamma}\\
			\qquad+\left\langle T\left( \begin{array}{ccc}
				\phi_j\\
				0
			\end{array}\right), \left( \begin{array}{ccc}
				\phi_i\\
				0
			\end{array}\right)\right\rangle_{\Gamma_b} ,  &1\leq i,j\leq m,\\
			\frac{\gamma}{\kappa^2}\langle \nabla \phi_{j-m},\nabla^\bot \phi_i\rangle_{\Omega_b}+\left\langle T\left( \begin{array}{ccc}
				0\\
				\phi_{j-m}
			\end{array}\right), \left( \begin{array}{ccc}
				\phi_i\\
				0
			\end{array}\right)\right\rangle_{\Gamma_b} ,  &1\leq i\leq m,m+1\leq j\leq 2m,\\
			\frac{\gamma}{\kappa^2}\langle \nabla \phi_{j},\nabla^\bot \phi_{i-m}\rangle_{\Omega_b}+\left\langle T\left( \begin{array}{ccc}
				\phi_{j}\\
				0
			\end{array}\right), \left( \begin{array}{ccc}
				0\\
				\phi_{i-m}
			\end{array}\right)\right\rangle_{\Gamma_b},  & m+1\leq i\leq 2m,1\leq j\leq m,\\
			\frac{\omega\mu}{\kappa^2}\langle \nabla \phi_{j-m},\nabla \phi_{i-m}\rangle_{\Omega_b}-\omega\mu \langle  \phi_{j-m}, \phi_{i-m}\rangle_{\Omega_b}\\
			\qquad+\frac{i}{\lambda}\langle  \phi_{j-m}, \phi_{i-m}\rangle_{\Gamma}+\left\langle T\left( \begin{array}{ccc}
				0\\
				\phi_{j-m}
			\end{array}\right), \left( \begin{array}{ccc}
				0\\
				\phi_{i-m}
			\end{array}\right)\right\rangle_{\Gamma_b},  & m+1\leq i,j\leq 2m,
		\end{cases}
	\end{align*}
	and let $ F\in \mathbb{C}^{2m} $ be a vector whose components are given by 
	\begin{align*}
		F_i=\begin{cases}
			-\frac{2i\omega\eps\beta e^{-i\beta b}}{\kappa^2}\int_{\Gamma_b}\eps p_3\overline{\phi}_i e^{i\alpha x_1}\,ds, & 1\leq i\leq m\\
			-\frac{2i\omega\eps\beta e^{-i\beta b}}{\kappa^2}\int_{\Gamma_b}\mu q_3\overline{\phi}_{i-m}e^{i\alpha x_1}\,ds, & m+1\leq i\leq 2m.
		\end{cases}
	\end{align*}
	Then we get the system of linear equations
	\begin{align}
		\sum_{j=1}^{2m}B_{ij}a_j=F_i, \quad 1\leq i\leq 2m.  \label{stiffness}
	\end{align}
	Having obtained $ \{a_j\}_{j=1}^{2m} $ from $ (\ref{stiffness}) $, we can get $ u_h $ and $ v_h $ by setting $ c_j=a_j,\;d_{j}=a_{j+m} $ for $ 1\leq j\leq m $.
	
	Below we prove the well-posedness of the finite element approximation problem $ (\ref{VFFEM}) $ and an error estimate of the finite element solution. 
	Denote $ e_h=(u-u_h,v-v_h) $. 
	It is obvious that $ e_h $ is $ \alpha $-quasiperiodic.
	Define the projection operator $ P:  L^2(\Gamma_b)^2 \to  L^2(\Gamma_b)^2$ by $$ (Pf)(x_1)= \sum_{n\in\mathcal{A}}f_n e^{i\alpha_n x_1}, \quad f=\sum_{n\in\mathbb{Z}}f_n e^{i\alpha_n x_1} \in L^2(\Gamma_b)^2,$$
	where the set $ \mathcal{A} $ is defined by $ (\ref{aA}) $.
	\begin{lemma}\label{FemH_1}
		There exists a constant $ h_1 \in (0,1) $ such that for $ h\in (0,h_1) $ the following estimate holds:
		\ben
		\|e_h\|_{H^1(\Omega_b)^2}^2 \leq C\left(h^{2\rho-2}\|(u,v)\|^2_{H^\rho(\Omega_b)^2}+\|e_h\|_{L^2(\Omega_b)^2}^2+\|Pe_h\|_{L^2(\Gamma_b)^2}^2\right),
		\enn
		where the constant $ C $ depends on $ \rho $ but is independent of $ h $ and $ (u,v) $.
	\end{lemma}
	\begin{proof}
		It follows from the sesquilinear form $ (\ref{VFFEM}) $ that
		\be
		B(e_h;e_h)&:=&
		\int_{\Omega_b}\left[\frac{\omega\eps}{\kappa^2}|\nabla (u-u_h)|^2-\frac{\gamma}{\kappa^2}\nabla (v-v_h)\cdot \nabla^\bot\overline{(u-u_h)}-\omega \eps\,|u-u_h|^2 \right.\nonumber\\
		&\,&+\left.\frac{\omega\mu}{\kappa^2}|\nabla (v-v_h)|^2
		+\frac{\gamma}{\kappa^2}\nabla (u-u_h)\cdot \nabla^\bot\overline{(v-v_h)}-\omega \mu\,|v-v_h|^2\right]\;dx \nonumber\\
		&\,&+\int_{\Gamma}\frac{i}{\lambda}|u-u_h|^2+i\lambda |v-v_h|^2\,ds+\int_{\Gamma_b}
		Te_h\cdot \overline{e_h}\,ds. \label{Be_h^2}
		\en
		Multiplying both sides of $ (\ref{Be_h^2}) $ by $ \theta=\frac{i+\delta}{|i+\delta|} $ and taking the real part, we get
		\be
		\mbox{Re}\,[\theta B(e_h;e_h)]&=&\mbox{Re}\,\left\{\theta\int_{\Omega_b}\left[\frac{\omega\eps}{\kappa^2}|\nabla (u-u_h)|^2-\frac{\gamma}{\kappa^2}\nabla (v-v_h)\cdot \nabla^\bot\overline{(u-u_h)}-\omega \eps\,|u-u_h|^2 \right.\right.\nonumber\\
		&\,&+\left.\left.\frac{\omega\mu}{\kappa^2}|\nabla (v-v_h)|^2
		+\frac{\gamma}{\kappa^2}\nabla (u-u_h)\cdot \nabla^\bot\overline{(v-v_h)}-\omega \mu\,|v-v_h|^2\right]\;dx \right\} \nonumber\\
		&\,&+\mbox{Re}\,\left\{\theta\int_{\Gamma}\frac{i}{\lambda}|u-u_h|^2+i\lambda |v-v_h|^2\,ds\right\}+\left\{\mbox{Re}\,\theta\int_{\Gamma_b}
		Te_h\cdot \overline{e_h}\,ds\right\}.\label{RealBe_h^2}
		\en
		From the strongly elliptic analysis (see the proof of Theorem \ref{StronglyEllip}), we can easily get 
		\ben
		\mbox{Re}\,\left\{\theta\int_{\Gamma}\frac{i}{\lambda}|u-u_h|^2+i\lambda |v-v_h|^2\,ds \right\}\geq 0,
		\enn
		and
		\ben
		\mbox{Re}\,\left\{\theta\int_{\Gamma_b}
		Te_h\cdot \overline{e_h}\,ds+q(e_h;e_h)\right\} \geq 0,
		\enn
		where the compact form $ q  $ is defined by (\ref{defQ}), that is, for $ e_h=\sum_{n\in\mathbb{Z}} A_n e^{i\alpha_n x_1} $, we have
		\ben
		q(e_h;e_h)=2\pi\mbox{Re}\,\sum_{n\in \mathcal{A}}\theta M_n A_n \cdot \overline{A_n} \leq C\|Pe_h\|_{L^2(\Gamma_b)^2}^2. 
		\enn
		Therefore, by (\ref{RealBe_h^2}),
		\ben
		&\quad&\mbox{Re}\,\left\{\theta\int_{\Omega_b}\left[\frac{\omega\eps}{\kappa^2}|\nabla (u-u_h)|^2-\frac{\gamma}{\kappa^2}\nabla (v-v_h)\cdot \nabla^\bot\overline{(u-u_h)}\right.\right. \nonumber\\
		&\;&+\left.\left.\frac{\omega\mu}{\kappa^2}|\nabla (v-v_h)|^2
		+\frac{\gamma}{\kappa^2}\nabla (u-u_h)\cdot \nabla^\bot\overline{(v-v_h)}\right]\;dx \right\}\nonumber\\
		&=&\mbox{Re}\,(\theta B(e_h;e_h))-\mbox{Re}\,\left\{\theta\int_{\Gamma}\frac{i}{\lambda}|u-u_h|^2+i\lambda |v-v_h|^2\,ds\right\}-\mbox{Re}\,\left\{\theta\int_{\Gamma_b}
		Te_h\cdot \overline{e_h}\,ds \right\}\nonumber\\
		&\;&+\mbox{Re}\,\left\{\theta\int_{\Omega_b}\omega \eps\,|u-u_h|^2+\omega \mu\,|v-v_h|^2 \;dx \right\}\nonumber\\
		&\leq& \mbox{Re}\,(\theta B(e_h;e_h))+\mbox{Re}\,\left\{\theta\int_{\Omega_b}\omega \eps\,|u-u_h|^2+\omega \mu\,|v-v_h|^2\;dx\right\}+q(e_h;e_h),
		\enn
		Using Lemma $ \ref{PD} $ (i) we get
		\begin{align}
			C_1\|e_h\|^2_{H^1(\Omega_b)^2} \leq |B(e_h;e_h)|+C_2\|e_h\|^2_{L^2(\Omega_b)^2}+C\|Pe_h\|_{L^2(\Gamma_b)^2}^2. \label{FemEsti11}
		\end{align}  
		Observing for any $ (\xi,\eta)\in X_h^2 $ that
		\ben
		B(u,v;\xi-u_h,\eta-v_h)=F(\xi-u_h,\eta-v_h),\quad B(u_h,v_h;\xi-u_h,\eta-v_h)=F(\xi-u_h,\eta-v_h),
		\enn
		we obtain 
		\begin{align}
			B(u-u_h,v-v_h;\xi-u_h,\eta-v_h)=0. \label{othgnal}
		\end{align}
		Therefore for any $\left(\xi, \eta\right) \in  X_h^2$, we have
		\begin{align}
			B(u-u_h,v-v_h;u-u_h,v-v_h)=B(u-u_h,v-v_h;u-\xi,v-\eta).  \label{femeqn}	
		\end{align}
		Since $ X_h^2 $ is of finite dimensions, it is complete and therefore closed. 
		Hence, the infimum in $ (\ref{asuptnXH}) $ is actually attained for $ (\varphi,\psi)=(u,v) $ in $ (\ref{asuptnXH}) $. 
		For any small positive constants $ \epsilon_i \,(i=1,2,3,4)\, $, it follows from (\ref{femeqn}) and Young's inequality that
		\be
		|B(e_h;e_h)|&=&|B(u-u_h,v-v_h;u-\xi,v-\eta)| \nonumber \\
		&=&\left|\int_{\Gamma}\frac{i}{\lambda}(u-u_h)\,\overline{(u-\xi)}+i\lambda (v-v_h)\,\overline{(v-\eta)}\,ds+\int_{\Omega_b}\frac{\omega\eps}{\kappa^2}\nabla (u-u_h)\cdot \nabla\overline{(u-\xi)}\right.\nonumber\\
		&\,&-\frac{\gamma}{\kappa^2}\nabla (v-v_h)\cdot \nabla^\bot\overline{(v-\eta)}-\omega \eps\,(u-u_h)\,\overline{(u-\xi)}+\frac{\omega\mu}{\kappa^2}\nabla (v-v_h)\cdot \nabla\overline{(v-\eta)}\nonumber\\
		&\,&+\frac{\gamma}{\kappa^2}\nabla (u-u_h)\cdot \nabla^\bot\overline{(v-\eta)}-\omega \mu\,(v-v_h)\,\overline{(v-\eta)}\;dx \nonumber\\
		&\,&+\left.\int_{\Gamma_b}
		T\left( \begin{array}{ccc}
			u-u_h\\
			v-v_h
		\end{array}\right)
		\cdot
		\left( \begin{array}{ccc}
			\overline{u-\xi}\\
			\overline{v-\eta}
		\end{array}\right)\,ds \right|\nonumber\\
		&\leq& \frac{1}{|\lambda|}\left(h\|u-u_h\|^2_{L^2(\Gamma)}+\frac{1}{h}\|u-\xi\|^2_{L^2(\Gamma)}\right)+|\lambda|\left(h\|v-v_h\|^2_{L^2(\Gamma)}+\frac{1}{h}\|v-\eta\|^2_{L^2(\Gamma)}\right)\nonumber\\
		&\,& +\frac{\omega\eps}{\kappa^2}\left(\epsilon_1 \|\nabla u-\nabla u_h\|^2_{L^2(\Omega_b)}+\frac{1}{4\epsilon_1}\|\nabla u-\nabla\xi\|^2_{L^2(\Omega_b)}\right)\nonumber\\
		&\,&+\frac{|\gamma|}{\kappa^2}\left(\epsilon_2 \|\nabla v-\nabla v_h\|^2_{L^2(\Omega_b)}+\frac{1}{4\epsilon_2}\|\nabla u-\nabla\xi\|^2_{L^2(\Omega_b)}\right)\nonumber\\
		&\,& +\omega \eps\left(h^2\|u-u_h\|^2_{L^2(\Omega_b)}+h^{-2}\|u-\xi\|^2_{L^2(\Omega_b)}\right)\nonumber\\
		&\,&+\frac{\omega\mu}{\kappa^2}\left(\epsilon_3 \|\nabla v-\nabla v_h\|^2_{L^2(\Omega_b)}+\frac{1}{4\epsilon_3}\|\nabla v-\nabla\eta\|^2_{L^2(\Omega_b)}\right)\nonumber\\
		&\,&+\frac{|\gamma|}{\kappa^2}\left(\epsilon_4 \|\nabla u-\nabla u_h\|^2_{L^2(\Omega_b)}+\frac{1}{4\epsilon_4}\|\nabla v-\nabla\eta\|^2_{L^2(\Omega_b)}\right)\nonumber\\
		&\,&+\omega \mu\left(h^2\|v-v_h\|^2_{L^2(\Omega_b)}+h^{-2}\|v-\eta\|^2_{L^2(\Omega_b)}\right)\nonumber\\
		&\,&+\left|\int_{\Gamma_b}
		T\left( \begin{array}{ccc}
			u-u_h\\
			v-v_h
		\end{array}\right)
		\cdot
		\left( \begin{array}{ccc}
			\overline{u-\xi}\\
			\overline{v-\eta}
		\end{array}\right)\,ds \right|.\label{FemEstiB}
		\en
		Using the continuity of the DtN map $ T $ (see Lemma \ref{ctnT}), trace theorem and Young's inequality, we have
		\be
		&\;&\left|\int_{\Gamma_b}
		T\left( \begin{array}{ccc}
			u-u_h\\
			v-v_h
		\end{array}\right)
		\cdot
		\left( \begin{array}{ccc}
			\overline{u-\xi}\\
			\overline{v-\eta}
		\end{array}\right)\,ds \right|\nonumber\\
		&\leq& \left \|T\left( \begin{array}{ccc}
			u-u_h\\
			v-v_h
		\end{array}\right)\right \|_{H^{-1/2}(\Gamma_b)^2} \left\|\left( \begin{array}{ccc}
			u-\xi\\
			v-\eta
		\end{array}\right)\right \|_{H^{1/2}(\Gamma_b)^2}\nonumber\\
		&\leq& C \|e_h\|_{H^{1/2}(\Gamma_b)^2}\left\|\left( \begin{array}{ccc}
			u-\xi\\
			v-\eta
		\end{array}\right)\right \|_{H^{1/2}(\Gamma_b)^2}\nonumber\\
		&\leq& C\left(\epsilon_5\|e_h\|^2_{H^1(\Omega_b)^2}+\frac{1}{4\epsilon_5}\left\|\left( \begin{array}{ccc}
			u-\xi\\
			v-\eta
		\end{array}\right)\right \|_{H^{1/2}(\Gamma_b)^2}\right).\label{FemEstiT}
		\en	
		One deduces from $ (\ref{asuptnXH}) $ and $ (\ref{FemEstiB}) $ - $ (\ref{FemEstiT}) $ that
		\begin{align}
			&\quad |B(e_h;e_h)|\nonumber\\
			&\leq Ch \|e_h\|^2_{L^2(\Gamma)^2}+\sigma \|e_h\|^2_{H^1(\Omega_b)^2}+C_1 h^2 \|e_h\|^2_{L^2(\Omega_b)^2}+C(\sigma) h^{2\rho-2}\|(u,v)\|^2_{H^\rho(\Omega_b)^2},\label{FemEstiBEnd}
		\end{align}
		where $ \sigma=\sigma(\epsilon_1,\cdots,\epsilon_5)>0 $.
		Combining $ (\ref{FemEsti11}) $ and $ (\ref{FemEstiBEnd}) $ leads to
		\ben
		C_1\|e_h\|^2_{H^1(\Omega_b)^2} &\leq &|B(e_h;e_h)|+C_2\|e_h\|^2_{L^2(\Omega_b)^2}+|q(e_h;e_h)| \nonumber\\
		&\leq& C_3 h \|e_h\|^2_{L^2(\Gamma)^2}+\sigma \|e_h\|^2_{H^1(\Omega_b)^2}+C_4 h^2 \|e_h\|^2_{L^2(\Omega_b)^2}\nonumber\\
		&\;&+C(\sigma) h^{2\rho-2}\|(u,v)\|^2_{H^\rho(\Omega_b)^2}
		+C_2\|e_h\|^2_{L^2(\Omega_b)^2}+C_5  \|Pe_h\|^2_{L^2(\Gamma_b)^2}. \label{eH1}
		\enn
		Using the estimate 
		\ben
		\|e_h\|_{L^2(\Gamma)^2}\leq  \|e_h\|_{H^{1/2}(\Gamma)^2} \leq C\|e_h\|_{H^1(\Omega_b)^2},\quad C>0,
		\enn	
		we get from (\ref{eH1}) that 
		\begin{align*}
			C_1\|e_h\|^2_{H^1(\Omega_b)^2} &\leq (C C_3  h + \sigma+C_4 h^2 )\|e_h\|^2_{H^1(\Omega_b)^2}+C(\sigma) h^{2\rho-2}\|(u,v)\|^2_{H^\rho(\Omega_b)^2}\\
			&\quad +C_2\|e_h\|^2_{L^2(\Omega_b)^2}+C_5  \|Pe_h\|^2_{L^2(\Gamma_b)^2}.
		\end{align*}
		Now	choose $ \sigma $ sufficienlty small and let $ h_1 $ be a constant such that $ \sigma+C C_3  h_1 +C_4 h_1^2< C_1 $. Then for $ h\in (0,h_1) $, we obtain the desired estimate of this lemma.

	\end{proof}
	We next estimate the $ L^2 $-norm of $ e_h $ in $ \Omega_b $.
	\begin{lemma}\label{FemL_2}
		There exists a constant $ h_2 \in (0,1)$ such that 
		\ben
		\|e_h\|_{L^2(\Omega_b)^2}\leq C(h+C_1 h^{3/2})\|e_h\|_{H^1(\Omega_b)^2} \quad \text{ for all }h\in (0,h_2),
		\enn
		where the constants $ C,\; C_1 $ depend on $ \rho $ but are independent of $ h $ and $ (u,v) $.
	\end{lemma}
	\begin{proof}
		We use the duality argument. 
		By the definition
		\be
		\|e_h\|_{L^2(\Omega_b)^2}=\sup_{(\phi,\zeta)\in C_0^{\infty}(\Omega_b)^2}\frac{(e_h,(\phi,\zeta))_{L^2(\Omega_b)^2}}{\|(\phi,\zeta)\|_{L^2(\Omega_b)^2}},\label{FemDualDef}
		\en
		where 
		\begin{align}
			(e_h,(\phi,\zeta))_{L^2(\Omega_b)^2}:=\int_{\Omega_b}\frac{\omega\eps}{\kappa^2}(u-u_h)\bar{\phi} +\frac{\omega\mu}{\kappa^2} (v-v_h)\bar{\zeta}\;dx.  \label{eHIner}
		\end{align}
		Consider a quasi-periodic solution $ (w,z) $ of the following boundary value problem:
		\begin{align}  
			\begin{array}{l}
				\qquad\left\{\begin{array}{ll}
					\Delta w+k^2w = -\bar{\phi}   &\mbox{in} \quad \Omega_b, \\
					\Delta z+k^2z = -\bar{\zeta}    &\mbox{in} \quad \Omega_b, \\
					\lambda\partial_n w-\frac{i\kappa^2}{\omega\epsilon}w+\frac{\lambda \gamma}{\omega\epsilon}\partial_\tau z =0       &\mbox{on} \quad \Gamma,\\
					\partial_n z-\frac{i\lambda\kappa^2}{\omega\mu}z-\frac{\lambda }{\omega\mu}\partial_\tau w =0       &\mbox{on} \quad \Gamma,\\
					T^*
					\left(\begin{array}{ccc}
						w\\
						z
					\end{array}\right)=\sum_{n\in\mathbb{Z}}M_n^*
					\left(\begin{array}{ccc}
						\hat{w}_n\\
						\hat{z}_n
					\end{array}\right) e^{i\alpha_n x_1 } & \mbox{on} \quad \Gamma_b,
				\end{array}\right. \\ \label{DualPlm}
			\end{array} 
		\end{align}
		where $ T^*$ is the adjoint operator of $ T $. 
		We can easily get the variational formulation of $ (\ref{DualPlm}) $ that for all $ (\varphi,\psi) \in X $,
		\begin{align}
			&\quad	\int_{\Omega_b}\left[\frac{\omega\eps}{\kappa^2}\nabla w\cdot \nabla\varphi-\frac{\gamma}{\kappa^2}\nabla z\cdot \nabla^\bot\varphi-\omega \eps\,w\;\varphi+\frac{\omega\mu}{\kappa^2}\nabla z\cdot \nabla\psi\right.
			+\left.\frac{\gamma}{\kappa^2}\nabla w\cdot \nabla^\bot\psi-\omega \mu\,z\;\psi\right]\;dx\nonumber\\
			&+\int_{\Gamma}\frac{i}{\lambda}w\,\varphi+i\lambda z\,\psi\,ds
			-\int_{\Gamma_b}T^*
			\left(\begin{array}{ccc}   
				w\\  
				z 
			\end{array}
			\right)
			\cdot
			\left( \begin{array}{ccc}
				\varphi\\
				\psi
			\end{array}\right)\,ds
			=\int_{\Omega_b}\frac{\omega\eps}{\kappa^2}\bar{\phi}\varphi+\frac{\omega\mu}{\kappa^2}\bar{\zeta}\psi. \label{VFdual}
		\end{align}
		Taking $ \varphi=u-u_h,\,\psi=v-v_h $ in (\ref{VFdual}), using the definition of $ (e_h,(\phi,\zeta))_{L^2(\Omega_b)^2} $ in (\ref{eHIner}) and recalling the form $ B(u,v;\varphi,\psi) $ in (\ref{VFB}), we get
		\begin{align}
			B(u-u_h,v-v_h;\bar{w},\bar{z})=(e_h,(\phi,\zeta))_{L^2(\Omega_b)^2}. \label{dualeqn1}
		\end{align}
		The well-posedness of the problem $ (\ref{DualPlm}) $ can be established by the same argument as the proof for the variational problem $ (\ref{VF}) $. 
		Moreover, we have
		\be
		\|(w,z)\|_{H^2(\Omega_b)^2}\leq C \|(\phi,\zeta)\|_{L^2(\Omega_b)^2}. \label{DualWposd}
		\en
		Using the orthogonal formula (\ref{othgnal}), we have
		\begin{align}
			B(u-u_h,v-v_h;\bar{w}-\xi,\bar{z}-\eta)
			&=B(u-u_h,v-v_h;\bar{w},\bar{z})-B(u-u_h,v-v_h;\xi,\eta)\nonumber\\
			&=B(u-u_h,v-v_h;\bar{w},\bar{z}).\label{dualeqn2}
		\end{align}
		Combining (\ref{dualeqn1}) and (\ref{dualeqn2}) gives  
		\be
		|(e_h,(\phi,\zeta))_{L^2(\Omega_b)^2}|=|B(e_h,(\bar{w},\bar{z}))|=|B(e_h,(\bar{w},\bar{z})-(\xi,\eta))|\quad \text{ for all } (\xi,\eta)\in X_h^2.\label{dualRelation}
		\en
		In particular, $ (\xi,\eta) $ can be chosen in such a way that the infimum is attained for $ (\varphi,\psi)=(w,z) $ in $ (\ref{asuptnXH}) $. 
		By arguing analogously to the proof of Lemma $ \ref{FemH_1} $, we deduce from $ (\ref{FemEstiB}) $ - $ (\ref{FemEstiT}) $ that
		\ben 
		&\;&|B(e_h,(\bar{w},\bar{z})-(\xi,\eta))|  \\
		&=& \left|\int_{\Gamma}\frac{i}{\lambda}(u-u_h)\,\overline{(\bar{w}-\xi)}+i\lambda (v-v_h)\,\overline{(\bar{z}-\eta)}\,ds+\int_{\Omega_b}\frac{\omega\eps}{\kappa^2}\nabla (u-u_h)\cdot \nabla\overline{(\bar{w}-\xi)}\right.\\
		&\,&-\frac{\gamma}{\kappa^2}\nabla (v-v_h)\cdot \nabla^\bot\overline{(\bar{w}-\xi)}-\omega \eps\,(u-u_h)\,\overline{(\bar{w}-\xi)}+\frac{\omega\mu}{\kappa^2}\nabla (v-v_h)\cdot \nabla\overline{(\bar{z}-\eta)}\\
		&\,&+\frac{\gamma}{\kappa^2}\nabla (u-u_h)\cdot \nabla^\bot\overline{(\bar{z}-\eta)}-\omega \mu\,(v-v_h)\,\overline{(\bar{z}-\eta)}\;dx \\
		&\,&+\left.\int_{\Gamma_b}
		T\left( \begin{array}{ccc}
			u-u_h\\
			v-v_h
		\end{array}\right)
		\cdot
		\left( \begin{array}{ccc}
			\overline{\bar{w}-\xi}\\
			\overline{\bar{z}-\eta}
		\end{array}\right)\,ds \right|.\\
		\enn 
		Using the Cauchy-Schwarz inequality, we continue to estimate the above equation by 
		\be
		&\;&|B(e_h,(\bar{w},\bar{z})-(\xi,\eta))| \nonumber \\
		&\leq& \frac{1}{|\lambda|}\left(h^{-1/2}\|u-u_h\|_{L^2(\Gamma)}h^{1/2}\|\bar{w}-\xi\|_{L^2(\Gamma)}\right)+|\lambda|\left(h^{-1/2}\|v-v_h\|_{L^2(\Gamma)}h^{1/2}\|\bar{z}-\eta\|_{L^2(\Gamma)}\right)\nonumber\\
		&\,& +\frac{\omega\eps}{\kappa^2}\left(\frac{1}{h} \|\nabla u-\nabla u_h\|_{L^2(\Omega_b)}h|\bar{w}-\xi|_{H^1(\Omega_b)}\right)+\frac{|\gamma|}{\kappa^2}\left(\frac{1}{h} \|\nabla v-\nabla v_h\|_{L^2(\Omega_b)}h|\bar{w}-\xi|_{H^1(\Omega_b)}\right)\nonumber\\
		&\,& +\omega \eps\left(\|u-u_h\|_{L^2(\Omega_b)}\|\bar{w}-\xi\|_{L^2(\Omega_b)}\right)+\frac{\omega\mu}{\kappa^2}\left(\frac{1}{h} \|\nabla v-\nabla v_h\|_{L^2(\Omega_b)}h|\bar{z}-\eta|_{H^1(\Omega_b)}\right)\nonumber\\
		&\,&+\frac{|\gamma|}{\kappa^2}\left(\frac{1}{h} \|\nabla u-\nabla u_h\|_{L^2(\Omega_b)}h|\bar{z}-\eta|_{H^1(\Omega_b)}\right)
		+\omega \mu\left(\|v-v_h\|_{L^2(\Omega_b)}\|\bar{z}-\eta\|_{L^2(\Omega_b)}\right)\nonumber\\
		&\,& +\frac{1}{h}\|e_h\|_{H^1(\Omega_b)^2}h\left\|\left( \begin{array}{ccc}
			u-\xi\\
			v-\eta
		\end{array}\right)\right \|_{H^{1/2}(\Gamma_b)^2}\nonumber\\
		&\leq&
		C_1 h\|(w,z)\|_{H^2(\Omega_b)^2}\left[h^{1/2}\|e_h\|_{H^1(\Omega_b)^2}+C_2\|\nabla e_h\|_{L^2(\Omega_b)^2}\right.\nonumber\\
		&\,&\left.+C_3h\|e_h\|_{L^2(\Omega_b)^2}+C_4\|e_h\|_{H^1(\Omega_b)^2}\right],\label{Numerator}
		\en
		where the constants $ C_j\,(\,j=1,2,3,4\,) $ depend on $ k $ but independent of $ h $. 
		Combining $ (\ref{FemDualDef}) $ and $ (\ref{DualWposd}),\,(\ref{dualRelation}) $ and $ (\ref{Numerator}) $, we can find a positive constant $ h_2\leq 1 $ such that for all $ h\in (0,h_2) $,
		\begin{align*}
			&\quad
			\|e_h\|_{L^2(\Omega_b)^2}=\sup_{(\phi,\zeta)\in C_0^{\infty}(\Omega_b)^2}\frac{(e_h,(\phi,\zeta))_{L^2(\Omega_b)^2}}{\|(\phi,\zeta)\|_{L^2(\Omega_b)^2}}\\
			&\leq\sup_{(\phi,\zeta)\in C_0^{\infty}(\Omega_b)^2}\frac{C_1 h\|(w,s)\|_{H^2(\Omega_b)^2}\left[h^{1/2}\|e_h\|_{H^1(\Omega_b)^2}+C_2\| e_h\|_{H^1(\Omega_b)^2}
				+C_3h\|e_h\|_{L^2(\Omega_b)^2}\right]}{\|(\phi,\zeta)\|_{L^2(\Omega_b)^2}} \\
			&\leq\sup_{(\phi,\zeta)\in C_0^{\infty}(\Omega_b)^2}\frac{C_1 h\|(\phi,\zeta)\|_{L^2(\Omega_b)^2}\left[h^{1/2}\|e_h\|_{H^1(\Omega_b)^2}+C_2\| e_h\|_{H^1(\Omega_b)^2}
				+C_3h\|e_h\|_{L^2(\Omega_b)^2}\right]}{\|(\phi,\zeta)\|_{L^2(\Omega_b)^2}} \\
			&= C h\|e_h\|_{H^1(\Omega_b)^2}+C_1 h^2\|e_h\|_{L^2(\Omega_b)^2}+C_2 h^{3/2}\|e_h\|_{H^1(\Omega_b)^2}.
		\end{align*}
		Now let $ h_2 $ be a constant that satisfies $ 1-C_1 h_2^2 >0 $. 
		We then obtain the estimate of this lemma for all $ h\in (0,h_2) $.
	\end{proof}
	We proceed with the estimate of the $ L^2 $-norm of $ Pe_h $ on $ \Gamma_b $.
	\begin{lemma}\label{FemPeh} 
		There exists a constant $ C $ such that
		\ben
		\|Pe_h\|_{L^2(\Gamma_b)^2} \leq C \|e_h\|_{L^2(\Omega_b)^2},
		\enn
		where the positive constant $ C $ is independent of $ h $ and $ (u,v) $.
	\end{lemma}
	\begin{proof}
		Define 
		$$ D=\{(x_1,x_2)\in \mathbb{R}^2: 0< x_1< 2\pi,b-\epsilon< x_2< b\}, $$ 
		where the constant $ \epsilon >0 $ is chosen to satisfy $ b-\epsilon >\max_{x\in \Gamma} \{x_2\} $.
		Suppose 
		$$ e_h(x) =\sum_{n\in\mathbb{Z}} A_n e^{i(\alpha_n x_1+\beta_n x_2)},\quad x \in D.$$
		Then $ P(e_h|_{\Gamma_b})=\sum_{n\in \mathcal{A}}A_n e^{i(\alpha_n x_1+\beta_n b)} $.
		Direct calculations show that 
		\ben
		\|e_h\|_{L^2(\Omega_b)^2}^2 &\geq& \|e_h\|_{L^2(D)^2}^2 =\int_{0}^{2\pi}\int_{b-\epsilon}^{b}|e_h|^2\;dx_2dx_1 \\
		&=&\int_{0}^{2\pi}\int_{b-\epsilon}^{b}\sum_{n\in\mathbb{Z}} A_n e^{i(\alpha_n x_1+\beta_n x_2)}\cdot \sum_{m\in\mathbb{Z}} \overline{A_m} e^{-i\alpha_m x_1-i\overline{\beta_m} x_2}\;dx_2dx_1 \\
		&=&\sum_{n,m\in\mathbb{Z}}A_n\cdot \overline{A_m}\int_{0}^{2\pi}e^{i(\alpha_n-\alpha_m) x_1}\;dx_1\int_{b-\epsilon}^{b}e^{i(\beta_n-\overline{\beta_m})x_2}\;dx_2\\
		&=& 2\pi \sum_{n\in\mathbb{Z}}|A_n|^2\int_{b-\epsilon}^{b}e^{i(\beta_n-\overline{\beta_n})x_2}\;dx_2\\
		&=& 2\pi \left( \sum_{|\alpha_n|\leq \kappa}\epsilon|A_n|^2+ \sum_{|\alpha_n|> \kappa}|A_n|^2\int_{b-\epsilon}^{b}e^{-2|\beta_n|x_2}\;dx_2\right)\\
		&=&2\pi \left( \sum_{|\alpha_n|\leq \kappa}\epsilon|A_n|^2+ \sum_{|\alpha_n|> \kappa}|A_n|^2 C_n\right),
		\enn
		where $ C_n=-\frac{\epsilon}{2|\beta_n|}\left(e^{-2|\beta_n|b}-e^{-2|\beta_n|(b-1)}\right)>0 $.
		By (\ref{aA}) and the proof of Theorem \ref{PD}, we can easily get that $ \mathcal{A} $ coincides with the set $ \{n\in\mathbb{Z}:\kappa \leq|\alpha_n|<k\} $.
		Then we have 
		\begin{align*}
			\mathcal{A}  \subset\mathcal{C}:=\{n\in \mathbb{Z}:|\alpha_n|\geq \kappa\}.
		\end{align*}
		Therefore,
		\begin{align*}
			2\pi \left( \sum_{|\alpha_n|\leq \kappa}\epsilon|A_n|^2+ \sum_{|\alpha_n|> \kappa}|A_n|^2 C_n\right)
			\geq 2\pi C\left(\sum_{n\in \mathcal{A}}|A_n|^2\right)=C	\|Pe_h\|^2_{L^2(\Gamma_b)^2} 
		\end{align*}
		where $ C=\min\{\epsilon,\min_{n\in \mathcal{A}}C_n\} $.
		
	\end{proof}
	The main result of this section is stated below.
	\begin{theorem}
		Suppose that $ (u,v) \in H^\rho(\Omega_b)^2,\, \rho\geq 2$, satisfies the variational problem $ (\ref{VF}) $. 
		Suppose also that the family of finite element spaces $ \{X_h^2\} $ satisfies the assumption $ (\ref{asuptnXH}) $. 
		Then there exists $ h_0 \in (0,1) $ such that for $ h \in (0, h_0) $, the problem $ (\ref{VFFEM}) $ admits a unique solution $ (u_h,v_h) $ with the estimates 
		\begin{align*}
			\|(u,v)-(u_h,v_h)\|_{L^2(\Omega_b)^2}&\leq C\left(h^\rho+C_1 h^{\rho+1/2}\right)\|(u,v)\|_{H^\rho(\Omega_b)^2},\\
			\|(u,v)-(u_h,v_h)\|_{H^1(\Omega_b)^2}&\leq C h^{\rho-1} \|(u,v)\|_{H^\rho(\Omega_b)^2},
		\end{align*}
		where the positive constant $ C $ depends on $ \rho $ but is independent of $ h $ and $ (u,v) $.
	\end{theorem}
	\begin{proof}
		Let $ h_1 $ and $  h_2 $ be specified as in Lemmas \ref{FemH_1} and \ref{FemL_2} and set $ h_0=\min \{h_1,h_2\} $. 
		For $ h\in (0,h_0) $, we deduce from Lemmas \ref{FemH_1}-\ref{FemPeh} that
		\begin{align*}
			\|e_h\|_{H^1(\Omega_b)^2}^2 
			&\leq C\left(h^{2\rho-2}\|(u,v)\|^2_{H^\rho(\Omega_b)^2}+C\|e_h\|^2_{L^2(\Omega_b)^2}\right)\\
			&\leq C\left(h^{2\rho-2}\|(u,v)\|^2_{H^\rho(\Omega_b)^2}+C_1(h+C_2 h^{3/2})^2\|e_h\|^2_{H^1(\Omega_b)^2}\right).
		\end{align*}
		Now	letting $ h_0 $ be a constant such that $ 1-C_1 (h_0+C_2 h_0^{3/2})^2>0  $, we obtain 
		\ben
		\|e_h\|_{H^1(\Omega_b)^2}\leq C h^{\rho-1} \|(u,v)\|_{H^\rho(\Omega_b)^2}\quad \text{ for all } h\in (0,h_0).
		\enn
		Therefore, using Lemma \ref{FemL_2},
		\ben
		\|e_h\|_{L^2(\Omega_b)^2}\leq C (h+C_1 h^{3/2})\|e_h\|_{H^1(\Omega_b)^2}\leq C\left(h^\rho+C_1 h^{\rho+1/2}\right)\|(u,v)\|_{H^\rho(\Omega_b)^2},
		\enn
		which completes the proof.
	\end{proof}
	
	\section{Integral equation methods} \label{sec:integralEm}
	The aim of this section is to develop an integral equation method for the conical diffraction problem (\ref{cdp}). 
	We make the following assumption. 
	\begin{description}
		\item \textbf{Assumption (A)}: The grating profile $\Gamma$ is the graph of some $2\pi$-periodic function $x_2=f(x_1)$, $x_1\in \R$, where $f$ is either $C^2$-smooth or piecewise linear with a finite number of corner points in one periodic cell.
	\end{description}
	Introduce the $\al$-quasiperiodic  fundamental solution to the Helmholtz equation $(\Delta+\kappa^2)u=0$ by
	\ben\no
	G(x,y)&=& \frac{i}{4} \sum_{n\in\Z} \exp(-i\al2\pi n)H^{(1)}_0\left(k\sqrt{(x_1+2n\pi-y_1)^2+(x_2-y_2)^2}\right)\\ 
	&=&\frac{i}{4\pi}\sum_{n\in\Z}\frac{1}{\beta_n}\exp\left(i \al_n(x_1-y_1)+i\beta_n|x_2-y_2|\right)
	\enn
	for $x-y\neq n(2\pi,0)$, with $H^{(1)}_0(t)$ being the first kind Hankel function of order zero. Define the single-layer potential by
	\ben
	(\mathcal{S} g)(x)= 2\int_{\Gamma}G(x,y)g(y)ds(y),\quad x\in \Om,
	\enn
	with the density $g$. We make the ansatz for the solution $(u^{s},v^{s})$ in the form
	\ben
	u^{s}=\mathcal{S}g_1,\quad  v^{s}=\mathcal{S}g_2.
	\enn
	Further, define the single- and double-layer operators $S$ and $K$ by
	\ben
	(S\rho)(x)&:=&2\int_\Gamma G(x,y)\rho(y)ds(y),\quad x\in\Gamma, \\
	(K\rho)(x)&:=&2\int_\Gamma\frac{\partial G(x,y)}{\partial \nu(y)}\rho(y)ds(y),\quad x\in\Gamma,
	\enn
	and the normal and tangential derivative operators $K'$ and $H'$ by
	\ben
	(K'\rho)(x)&:=&2\int_\Gamma\frac{\partial G(x,y)}{\partial \nu(x)}\rho(y)ds(y),\quad
	x\in\Gamma,\\
	(H'\rho)(x)&:=&2\int_\Gamma\frac{\partial G(x,y)}{\partial \tau(x)}\rho(y)ds(y),\quad
	x\in\Gamma,
	\enn
	where $ \nu  $ denotes the unit normal vector to the boundary $ \Gamma $ directed into the exterior of $ \Omega  $ and $ \tau $ denotes the unit tangential vector to $ \Gamma $.
	\begin{lemma}
		Let $ g_1,\, g_2 $ be the density functions of $u^{s},\, v^{s}$, respectively, then the following jump relations hold
		\begin{align*}
			u^s(x)&=2\int_\Gamma G(x,y)g_1(y)ds(y)=S g_1,\quad x\in\Gamma, \\
			v^s(x)&=2\int_\Gamma G(x,y)g_2(y)ds(y)=S g_2,\quad x\in\Gamma, \\
			\frac{\partial u^s_{\pm}}{\partial \nu}(x)&=2\int_\Gamma\frac{\partial G(x,y)}{\partial \nu(x)} g_1(y)ds(y) \pm  g_1(x)=K'g_1(x)\pm  g_1(x),\quad
			x\in\Gamma,\\
			\frac{\partial v^s_{\pm}}{\partial \nu}(x)&=2\int_\Gamma\frac{\partial G(x,y)}{\partial \nu(x)} g_2(y)ds(y) \pm  g_2(x)=K'g_2(x)\pm  g_2(x),\quad
			x\in\Gamma,\\
			\frac{\partial u^s}{\partial \tau}(x)&=2\int_\Gamma\frac{\partial G(x,y)}{\partial \tau(x)} g_1(y)ds(y)=H'g_1(x),\quad x\in\Gamma, \\
			\frac{\partial v^s}{\partial \tau}(x)&=2\int_\Gamma\frac{\partial G(x,y)}{\partial \tau(x)} g_2(y)ds(y)=H'g_2(x),\quad x\in\Gamma, 
		\end{align*}
		where
		\begin{align*}
			\frac{\partial u^s_{\pm}}{\partial \nu}(x)&:=\lim_{h\to +0} \nu(x)\cdot \nabla u^s(x\pm h\nu(x)), \quad
			\frac{\partial v^s_{\pm}}{\partial \nu}(x):=\lim_{h\to +0} \nu(x)\cdot \nabla v^s(x\pm h\nu(x)),\\
			\frac{\partial u^s}{\partial \tau}(x)&:=\tau (x)\cdot \nabla u^s(x), \qquad \qquad\qquad \quad\quad
			\frac{\partial v^s}{\partial \tau}(x):=\tau (x)\cdot \nabla v^s(x).
		\end{align*}
	\end{lemma}
	\begin{proof}
		We refer to \cite{MR4079757} for the mapping properties of the single- and double-layer operators.
		By reference \cite{Kress}, we only need to prove the continuity of $ \frac{\partial u}{\partial \tau}(x) $ on $ \Gamma $.
		\begin{align*}
			\lim_{x\to \Gamma^\pm} \nabla u^s(x)=2\int_\Gamma \nabla_x G(x,y)g_1(y)ds(y) \pm g_1(y)\nu(y).
		\end{align*}
		Therefore,
		\begin{align*}
			\lim_{x\to \Gamma^\pm}\frac{\partial u^s}{\partial \tau}(x)&=\lim_{x\to \Gamma^\pm}\tau (x)\cdot \nabla u^s(x)\\
			&=\tau (x) \cdot 2\int_\Gamma \nabla_x G(x,y)g_1(y)ds(y) \pm \tau (x)\cdot g_1(y)\nu(y)\\
			&=2\int_\Gamma\frac{\partial G(x,y)}{\partial \tau(x)} g_1(y)ds(y).
		\end{align*}
	\end{proof}
	Using the jump relations for (\ref{BC1}), we have on $ \Gamma $ that
	\begin{align}	
		&\quad\la\frac{\partial u^s}{\partial \nu}+i\omega\mu \cos^2\phi\; u^s+\la\sin\phi\sqrt{\frac{\mu}{\eps}}\frac{\partial v^s}{\partial\tau} \nonumber\\
		&=\la(K'g_1+g_1)+i\omega\mu \cos^2\phi Sg_1+\la\sin\phi\sqrt{\frac{\mu}{\eps}}H'g_2=h_1,\label{JR1}\\
		&\quad \frac{\partial v^s}{\partial n}+i\la\omega\eps\cos^2\phi\, v^s-\sin\phi\sqrt{\frac{\eps}{\mu}}\frac{\partial u^s}{\partial \tau}\nonumber\\
		&=(K'g_2+g_2)+i\la\omega\eps\cos^2\phi Sg_2 -\sin\phi\sqrt{\frac{\eps}{\mu}}H'g_1=h_2,\label{JR2}
	\end{align}
	where
	\begin{align*}
		h_1&:=-\left(\la\frac{\partial u^i}{\partial \nu}+i\omega\mu \cos^2\phi\; u^i+\la\sin\phi\sqrt{\frac{\mu}{\eps}}\frac{\partial v^i}{\partial\tau}\right),\\
		h_2&:=-\left(\frac{\partial v^i}{\partial \nu}+i\la\omega\eps\cos^2\phi\, v^i-\sin\phi\sqrt{\frac{\eps}{\mu}}\frac{\partial u^i}{\partial \tau}\right).
	\end{align*}
	Combining (\ref{JR1}) with (\ref{JR2}), we obtain the integral equations
	\ben
	\begin{pmatrix}
		\la(K'+I) & \la\sin\phi\sqrt{\mu/\eps} H'\\ 
		-\sin\phi\sqrt{\eps/\mu} H' & K'+I
	\end{pmatrix} \begin{pmatrix}
		g_1\\
		g_2
	\end{pmatrix}  + \begin{pmatrix}
		i \omega\mu\cos^2\phi S & 0 \\
		0 & i\la\omega\eps\cos^2\phi S \end{pmatrix} \begin{pmatrix}
		g_1\\
		g_2
	\end{pmatrix}=\begin{pmatrix}
		h_1\\
		h_2
	\end{pmatrix}.
	\enn 
	Therefore the equivalent system is
	\be\label{int}
	Ag+Bg:=\begin{pmatrix}
		\la(K'+I) & d H'\\ 
		-c H' & K'+I
	\end{pmatrix} g  + \begin{pmatrix}
		i a S & 0 \\
		0 & i b S \end{pmatrix} g=h, 
	\en with $g=(g_1,g_2)^\top, h=(h_1,h_2)^\top\in H^{-1/2}(\Gamma)^2$, and
	\ben
	d=\la\sin\phi\sqrt{\mu/\eps},\; c= \sin\phi\sqrt{\eps/\mu},\;a=\omega\mu\cos^2\phi,\;b=\la\omega\eps\cos^2\phi.
	\enn
	
	Under the Assumption (A), the single-layer operator $S$ is invertible form $H^{-1/2}(\Gamma)\rightarrow H^{1/2}(\Gamma)$.
	
	\begin{theorem}\label{ITH2} 
		Suppose that Assumption (A) holds. Then the operator $A+B$ defined in (\ref{int}) is a Fredholm operator with an index zero. 
		Moreover, the system (\ref{int}) admits a unique solution if $k^2\neq\gamma^2$ and $\lambda<0$. 
	\end{theorem}
	\begin{proof}
		It suffices to prove the Fredholm property of $A+B$, since the second assertion of Theorem \ref{ITH2} follows directly from the Fredholm alternative combined with Theorem
		\ref{TH1}. 
		To do this, we introduce the adjoint operator $H$ of $H'$, given by
		\ben 
		(Hg)(x)&:=&2\int_\Gamma\frac{\partial G(x,y)}{\partial \tau(y)}g(y)ds(y),\quad
		x\in\Gamma.
		\enn
		It is known that the adjoint operator of $K$ is just $K'$. 
		Since the operator $S$ is compact from $H^{-1/2}(\Gamma)\rightarrow H^{-1/2}(\Gamma)$, we only need to justify the Fredholm property of the adjoint operator $A^*$ of $A$, given by
		\ben
		A^*=\begin{pmatrix}
			\la(I+K) & -c H\\ 
			d H & I+K
		\end{pmatrix}:\quad H^{1/2}(\Gamma)^2\rightarrow H^{1/2}(\Gamma)^2.
		\enn
		It is easy to see that the operator $H_1=H+j$ with the rank $1$ operator 
		\ben
		j u=(u,e)_{L^2(\Gamma)}\; e,\quad e=1\in \C, 
		\enn
		is invertible in $H^{1/2}(\Gamma)$. 
		We will show that the operator
		\ben
		A_1:=\begin{pmatrix}
			\la(I+K) & -c H_1\\ 
			d H_1 & I+K
		\end{pmatrix}:\quad H^{1/2}(\Gamma)^2\rightarrow H^{1/2}(\Gamma)^2
		\enn 
		is a Fredholm operator with an index zero. 
		Simple calculations show that the operator 
		\ben
		B_1:=\begin{pmatrix}
			-(dH_1)^{-1}(I+K) & I \\ I & 0
		\end{pmatrix}=
		\begin{pmatrix}
			0 & I \\ I & (dH_1)^{-1}(I+K)
		\end{pmatrix}^{-1},
		\enn 
		is invertible, and that 
		\be\label{I0}
		A_1B_1=\begin{pmatrix}
			-\la(I+K)(dH_1)^{-1}(I+K)-cH_1& \la(I+K)\\ 0 & dH_1 
		\end{pmatrix}.
		\en
		Using the relations $HK=-KH$ and $(I+K)\,e=0$ (see \cite{Kress}), we get
		\ben
		(I+K)H_1&=&H_1(I-K)-j (I-K),\enn and thus
		\be \label{I1}
		(dH_1)^{-1}(I+K)&=&d^{-1}(I-K)H_1^{-1}-(dH_1)^{-1}[j(I-K)]H_1^{-1}.
		\en 
		Inserting (\ref{I1}) into (\ref{I0}) gives
		\ben 
		A_1B_1=\begin{pmatrix}
			-d^{-1}[\la(I-K^2)+cd H_1^2]H_1^{-1}+j_1& \la(I+K)\\ 0 & dH_1 
		\end{pmatrix},
		\enn 
		with $j_1:=\la (I+K)(dH_1)^{-1}[j(I-K)]H_1^{-1}$ being a rank one operator. 
		Hence, $A_1$ is Fredholm with index zero if this is true for the operator $\la(I-K^2)+cd H_1^2$.
		Making use of $K^2-H^2=I$ and the definitions of $c$ and $d$, we find
		\be\label{I5}
		\la(I-K^2)+cd H_1^2=(cd-\la)\, H_1^2+j_2=-\la\,\cos^2\phi\,  H_1^2+j_2,\en
		where $j_2$ is some operator with rank one. 
		Since $|\phi|<\pi/2$, we finally conclude that the operator (\ref{I5}) is  Fredholm  with index zero. 
		Theorem \ref{ITH2} is thus proven.
	\end{proof}
	
	\section*{Acknowledgments}
	
	\bibliographystyle{plain}
	\bibliography{refs}
\end{document}